\documentclass{m2an}
\newcommand{\gvec}{\ensuremath{\mathbf{g}}}
\newcommand{\mvec}{\ensuremath{\mathbf{m}}}
\newcommand{\nvec}{\ensuremath{\mathbf{n}}}
\newcommand{\uvec}{\ensuremath{\mathbf{u}}}
\newcommand{\vvec}{\ensuremath{\mathbf{v}}}
\newcommand{\xvec}{\ensuremath{\mathbf{x}}}
\newcommand{\yvec}{\ensuremath{\mathbf{y}}}

\newcommand{\Div}{\ensuremath{\mathrm{div}}}
\begin{document}
\title{Solvability of the mixed formulation for \\
       Darcy--Forchheimer flow in porous media}
\runningtitle{Mixed formulation for Darcy--Forchheimer flow}
\author{Peter Knabner}%
\address{Institute for Applied Mathematics,
	 Martensstra{\ss}e 3, D-91058 Erlangen, Germany;
         \email{knabner@am.uni-erlangen.de \&\ summ@am.uni-erlangen.de}}
\author{Gerhard Summ}\sameaddress{1}
\date{\today}
\begin{abstract} 
We consider the mixed formulation of the equations governing
Darcy--Forchheimer flow in porous media.
We prove existence and uniqueness of a solution for the stationary problem
and the existence of a solution for the transient problem.
 \end{abstract}
\begin{resume} Nous \'etudions la formulation mixte des \'equations 
pour l'\'ecoulement d'un gaz \`a travers d'un milieu poreux, 
qu'on suppose r\'egi par la loi de Darcy--Forchheimer.
Nous \'etablions des r\'esultats d'existence et d'unicit\'e pour le probl\`eme
stationnaire et un r\'esultat d'existence pour le probl\`eme transitoire.
\end{resume}
\maketitle
\section*{Introduction}
The flow of a gas through a porous medium is governed 
by the doubly nonlinear parabolic equation
\begin{equation}				\label{parabol.eq}
\phi(\xvec) \, \partial_t \rho(S(\xvec,t),\xvec,t) -
\Div \left( F (\nabla S(\xvec,t),\xvec,t)\right) = f(\xvec,t) \: ,
\quad (\xvec,t) \in \Omega \times [0,T] \: ,
\end{equation}
where the unknown function $S (= |p|p)$ represents the pressure squared
and the nonlinearities $\rho$ and $F$ are defined by
\begin{eqnarray}				\label{def:rho}
\rho(S(\xvec,t),\xvec,t) & := & 
\gamma(\xvec,t) \frac{S(\xvec,t)}{\sqrt{|S(\xvec,t)|}} \: , \\ \label{def:F}
F ( \nabla S(\xvec,t),\xvec,t ) & := & \frac{\sqrt{\alpha(\xvec,t)^2 + 
		 4 \beta(\xvec,t) |\nabla S(\xvec,t)|} - \alpha(\xvec,t)}
		{2 \beta(\xvec,t) |\nabla S(\xvec,t)|} \nabla S(\xvec,t) \: .
\end{eqnarray}
This equation together with appropriate initial and Neumann boundary conditions
has been studied by Amirat \cite{Amirat:91}.
He restricts his considerations to the case, where $\gamma$ is constant,
and shows the existence of a solution using the technique 
of semi-discretization in time.
Under additional regularity conditions on the solution he proves
the uniqueness and positivity of this solution.
The technique of semi-discretization in time has been used several times
to study similar doubly nonlinear parabolic equations.
We mention only the articles of Raviart, e.g.\ \cite{Raviart:70}.

Equation (\ref{parabol.eq}) can be derived from the following system of 
equations consisting of the Darcy--Forchheimer equation
(see e.g.\ \cite{Whitaker:96})
$$ 
\frac{\mu(\xvec,t)}{k(\xvec)} \, \uvec(\xvec,t)
+ \beta_{\mathrm{Fo}}(\xvec) \, \rho(\xvec,t) \, 
  |\uvec(\xvec,t)| \, \uvec(\xvec,t)
+ \nabla p(\xvec,t) = 0 \: ,
\quad (\xvec,t) \in \Omega \times [0,T] \: ,
$$ 
the continuity equation
$$ 
\phi(\xvec) \, \frac{\partial \rho(\xvec,t)}{\partial t} 
+ \Div(\rho(\xvec,t) \uvec(\xvec,t)) = f(\xvec,t) \: ,
\quad (\xvec,t) \in \Omega \times [0,T]
$$ 
and the ideal gas law as equation of state
$$ 
\rho(\xvec,t) = \frac{p(\xvec,t) \, W(\xvec,t)}{R_0 \, \Theta(\xvec,t)} 
              =: p(\xvec,t) \gamma(\xvec,t)
\: ,
\quad (\xvec,t) \in \Omega \times [0,T] \: .
$$ 
The unknowns here are the pressure $p$, the density $\rho$ 
and the volumetric flow rate $\uvec$ of the gas.
Porosity $\phi$, permeability $k$ 
and Forchheimer coefficient $\beta_{\mathrm{Fo}}$ of the porous medium,
viscosity $\mu$, molecular weight $W$ and temperature $\Theta$ of the gas,
and the universal gas constant $R_0$ are given as well as the source term $f$.
Assuming $\rho>0$ and introducing new variables $S = |p| p$ and 
$\mvec = |\rho| \uvec$ these equations can be transformed into
\begin{eqnarray}				\label{Da-Fo.eq}
\left( \alpha(\xvec,t) + \beta(\xvec,t) |\mvec(\xvec,t)| \right) \mvec(\xvec,t)
+ \nabla S(\xvec,t) & = & 0 \: ,
\hspace{3.25em} (\xvec,t) \in \Omega \times [0,T] \: , \\	\label{cont.eq}
\phi(\xvec) \, \partial_t \rho(S(\xvec,t),\xvec,t) 
+ \Div \left( \mvec(\xvec,t) \right) & = & f(\xvec,t) \: ,
\quad (\xvec,t) \in \Omega \times [0,T] \: ,
\end{eqnarray}
where 
$$ \gamma(\xvec,t) := \frac{W(\xvec,t)}{R_0 \, \Theta(\xvec,t)} \; , \quad
   \alpha(\xvec,t) := \frac{2 \, \mu(\xvec,t)}{\gamma(\xvec,t) \, k(\xvec)} \; , \quad
   \beta(\xvec,t) := \frac{2 \, \beta_{\mathrm{Fo}}(\xvec)}{\gamma(\xvec,t)} \; , $$
and the equation of state $\rho=\rho(S)$ is defined in (\ref{def:rho}).
Evidently, the Darcy--Forchheimer equation (\ref{Da-Fo.eq}) can be resolved
to give $\mvec = F (\nabla S)$ 
with the nonlinear mapping $F$ defined in (\ref{def:F}).
Substituting $\mvec = F (\nabla S)$ into (\ref{cont.eq}) finally yields
the parabolic equation (\ref{parabol.eq}).

The corresponding stationary problem of system (\ref{Da-Fo.eq}--\ref{cont.eq})
together with Neumann boundary conditions
has been studied by Fabrie \cite{Fabrie:89} for constant physical parameters.
He obtains the existence and uniqueness of a solution 
$(\mvec,S) \in \left( \xLn{3}(\Omega) \right)^n \times \mathrm{W}^{1,3/2}(\Omega)$ 
and shows additional regularity properties.

Throughout this article, for $s \in [0,\infty]$ we denote 
the Lebesgue spaces of $s$-integrable functions by $\xLn{s}(\Omega)$, 
for $m \ge 0, s \in [0,\infty]$ the Sobolev spaces by $\xWn{{m,s}}(\Omega)$
and the norm of $\xWn{{m,s}}(\Omega)$ by $\|\cdot\|_{m,s,\Omega}$ 
(\cf \cite{Adams}).
We assume that the domain $\Omega$ is bounded and fulfills
the uniform $\xCone$-regularity property.
Then the trace operator 
$\gamma_0: \xWn{{m,s}}(\Omega) \to \xWn{{m-1/s,s}}(\partial \Omega)$
is onto \cite[Thm.~7.53]{Adams}.
We denote by $\xWn{{m,s}}_0(\Omega)$ the kernel of $\gamma_0$ and its dual 
space by $\xWn{{-m,r}}(\Omega) := \left( \xWn{{m,s}}_0(\Omega) \right)'$,
where $1/s + 1/r = 1$.
Note that for every $s \in [0,\infty], m \ge 0$ the test space 
$\mathcal{D}(\Omega) := \xCinfty_0(\Omega)$ is a dense subset of 
$\xLn{s}(\Omega)$ and of $\xWn{{m,s}}_0(\Omega)$, and that 
$\mathcal{D}(\bar{\Omega}) := 
 \left\{ \Psi|_\Omega \bigm| \Psi \in \mathcal{D}(\xR^n) \right\}$
is a dense subset of $\xWn{{m,s}}(\Omega)$.
In addition, let us introduce the generalization $\xWn{s}(\Div;\Omega)$ of 
$\mathrm{H}(\Div;\Omega)$, which is defined in Appendix~\ref{app:W^s(div)}.
For $s=3$, this space will turn out to be appropriate for the nonlinearity
of the Darcy--Forchheimer law (see Proposition \ref{prop:A}).
Finally, we employ the spaces $\mathrm{C} \left( [0,T] ; X \right)$
and $\xLn{s} \left( 0,T; X \right)$ of vector-valued functions,
where $X$ is one of the above introduced spaces.

Note that the restriction of \cite{Fabrie:89} prevents from generalization 
to relevant situations from applications, where the parameters, 
e.g.\ the porosity $\phi$, vary discontinuously due to composite media
or where some of them, e.g. the temperature $\Theta$, are unknowns of 
a more complex model.
Both situations appear in the modelling of combustion in porous media,
see \cite{??}, which is our final goal.
Therefore we study the system (\ref{Da-Fo.eq}--\ref{cont.eq}) 
for Dirichlet boundary conditions and general coefficient functions,
imposing only minimal regularity assumptions.
We consider the stationary problem in Section \ref{sec:stat.pr} and
prove the existence and uniqueness of a solution.
To this end, we use a regularization of the equations and exploit
the monotonicity of the nonlinear mapping $F$.
In Section \ref{sec:semid.pr} we investigate the semi-discrete problem
after discretization of the time-derivative in (\ref{cont.eq})
and show again the existence and uniqueness of a solution.
Finally, in Section \ref{sec:trans.pr}, we study the transient problem
governed by (\ref{Da-Fo.eq}--\ref{cont.eq}) and prove its solvability.
Owing to the minor regularity of the solution of the transient problem,
we restrict our considerations there to homogeneous boundary conditions.

%
\section{The stationary problem}
\label{sec:stat.pr}
%
We consider the stationary problem governed by the Darcy--Forchheimer equation
and the stationary continuity equation together with Dirichlet boundary 
conditions:
\begin{eqnarray}				\nonumber
\left( \alpha(\xvec) + \beta(\xvec) |\mvec(\xvec)| \right) \mvec(\xvec)
+ \nabla S(\xvec) & = & 0 \: , 
\hspace{2.85em} \xvec \in \Omega \: , \\	\label{stat.pr}
\Div \left( \mvec(\xvec) \right) & = & f(\xvec) \: , 
\hspace{1.4em} \xvec \in \Omega \: , \\		\nonumber
S(\xvec) & = & S_b(\xvec) \: , \quad \xvec \in \partial \Omega \: .
\end{eqnarray}
We require $f \in \xLn{3}(\Omega)$, 
$S_b \in \mathrm{W}^{1/3,3/2}(\partial\Omega)$,
$\alpha, \beta \in \xLinfty(\Omega)$ and additionally 
$$ \left. \begin{array}{l}
0 < \underline{\alpha} \le \alpha(\xvec) \le \overline{\alpha} < \infty \: , \\
0 < \underline{\beta} \le \beta(\xvec) \le \overline{\beta} < \infty
\end{array} \right\}
\mbox{ for almost every } \xvec \in \Omega \; . $$

%
\subsection{Mixed formulation of the stationary problem}
The mixed formulation of (\ref{stat.pr}) reads as follows:
Find $(\mvec,S) \in \xWn{3}(\Div;\Omega) \times \xLn{{3/2}}(\Omega)$ such that
\begin{equation} 			\label{var.form.stat.pr}
\begin{array}{rcl@{\quad}l}
\displaystyle \int_\Omega \left( \alpha + \beta |\mvec| \right) 
               (\mvec \cdot \vvec) \, d\xvec
- \int_\Omega \Div(\vvec) \, S \, d\xvec & = &
\displaystyle - \int_{\partial\Omega} S_b (\vvec \cdot \nvec) \, d\sigma 
& \mbox{for all } \vvec \in \xWn{3}(\Div;\Omega) \: , \\[2ex]
\displaystyle \int_\Omega \Div(\mvec) \, q \, d\xvec & = & 
\displaystyle \int_\Omega f q \, d\xvec
& \mbox{for all } q \in \xLn{{3/2}}(\Omega) \: . 
\end{array}
\end{equation}
Next, we introduce continuous linear forms 
$g: \xWn{3}(\Div;\Omega) \to \xR$ and $f: \xLn{{3/2}}(\Omega) \to \xR$
by means of
$$ g(\vvec) := - \int_{\partial\Omega} S_b (\vvec \cdot \nvec) \, d\sigma \quad
   \mbox{for } \vvec \in \xWn{3}(\Div;\Omega) \quad , \quad
   f(q) := \int_\Omega f q \, d\xvec \quad
   \mbox{for } q \in \xLn{{3/2}}(\Omega) \: , $$
a bilinear form 
$b: \xWn{3}(\Div;\Omega) \times \xLn{{3/2}}(\Omega) \to \xR$
and a nonlinear form 
$a: \left( \xLn{3}(\Omega) \right)^n \times \left( \xLn{3}(\Omega) \right)^n 
    \to \xR$
by means of
$$ b(\vvec,q) := \! \int_\Omega \! \Div(\vvec) \, q \, d\xvec
   \mbox{ for } \vvec \in \xWn{3}(\Div;\Omega) \: , ~ q \in \xLn{{3/2}}(\Omega)
   ~ , ~
   a(\uvec,\vvec) := \! \int_\Omega \! \left( \alpha + \beta |\uvec| \right)
                                 (\uvec \cdot \vvec) \, d\xvec
   \mbox{ for } \uvec,\vvec \in \left( \xLn{3}(\Omega) \right)^n . $$
The form $b$ obviously is continuous, the continuity of the form $a$ will be
shown in Proposition \ref{prop:A}.
Then we can write the mixed formulation (\ref{var.form.stat.pr}) 
of (\ref{stat.pr}) in the following way:
Find $(\mvec,S) \in \xWn{3}(\Div;\Omega) \times \xLn{{3/2}}(\Omega)$ such that
\begin{equation} 			\label{mix.form.stat.pr}
\begin{array}{rcl@{\quad}l}
a(\mvec,\vvec) - b(\vvec,S) & = & g(\vvec) 
& \mbox{for all } \vvec \in \xWn{3}(\Div;\Omega) \: , \\[0.5ex]
b(\mvec,q) & = & f(q) 
& \mbox{for all } q \in \xLn{{3/2}}(\Omega) \: .
\end{array}
\end{equation}
In the following let $V$ be one of the spaces $\xWn{3}(\Div;\Omega)$ or
$\left( \xLn{3}(\Omega) \right)^n$.
Since $a$ is linear with respect to its second variable,
we can define a mapping $A: V \to V'$ by 
$\langle A\uvec , \vvec \rangle_{V' \times V} = a(\uvec,\vvec)$,
where $\langle \cdot , \cdot \rangle_{V' \times V}$ denotes the dual pairing
between $V'$ and $V$.
Using H\"older's inequality we obtain the following bound on $\|A\uvec\|_V$:
\begin{equation}			\label{norm(Au)}
\| A\uvec \|_{V'} \le C(\overline{\alpha}) \| \uvec \|_V + 
                         C(\overline{\beta}) \| \uvec\|_V^2 \: ,
\end{equation}
where $C(\overline{\alpha})$ and $C(\overline{\beta})$ are constants depending
only on $\overline{\alpha}$, $\overline{\beta}$ and the domain $\Omega$.
Thus $A\uvec$ is a continuous linear form on $V$ for every $\uvec \in V$.
The proof of (\ref{norm(Au)}) is contained in the proof of Proposition \ref{prop:A}.

The proof of the continuity and monotonicity of the mapping $A$ 
is based on the following lemma:
\begin{lmm}
The following inequalities hold for every $\xvec, \yvec \in \xR^n$:
\begin{eqnarray}				\label{ineq.1}
\left| |\xvec| \xvec - |\yvec| \yvec \right|
& \le & \left( |\xvec| + |\yvec| \right) \left| \xvec - \yvec \right| \: , \\
\left( |\xvec| \xvec - |\yvec| \yvec \right) \cdot \left( \xvec - \yvec \right)
& \ge & \frac{1}{2} |\xvec - \yvec|^3 \: .	\label{ineq.2}
\end{eqnarray}
\end{lmm}
The proof of (\ref{ineq.1}) and (\ref{ineq.2}) is elementary.
In the two-dimensional case, similar inequalities for more general 
nonlinearities are derived in \cite[Section~5]{Glowinski/Marroco:75}.

In the following we use the notions of \cite[Def.~25.2]{Zeidler}.
\begin{prpstn}				\label{prop:A}
The operator $A: V \to V'$ is continuous and strictly monotone on 
$V = \xWn{3}(\Div;\Omega)$, and continuous, uniformly monotone and coercive on 
$V = \left( \xLn{3}(\Omega) \right)^n$.
\end{prpstn}
\begin{proof}
To show the continuity of $A$ we consider 
$\uvec_1,\uvec_2 \in \left( \xLn{3}(\Omega) \right)^n$. 
Applying H\"older's inequality we obtain 
$$ \left| \langle A\uvec_1 - A\uvec_2 , \vvec \rangle_{V \times V'} \right| \le
   \left( C(\overline{\alpha}) \|\uvec_1-\uvec_2\|_{0,3,\Omega}
      + \overline{\beta} 
        \left\| |\uvec_1|\uvec_1-|\uvec_2|\uvec_2 \right\|_{0,3/2,\Omega}
          \right) \|\vvec\|_{0,3,\Omega} \quad
   \mbox{for all } \vvec \in \left( \xLn{3}(\Omega) \right)^n \: , $$
which also proves (\ref{norm(Au)}) for $V=\left( \xLn{3}(\Omega) \right)^n$
and therefore also for $V=\xWn{3}(\Div;\Omega)$.
Applying inequality (\ref{ineq.1}) and again the H\"older's inequality yields
$$ \left\| A \uvec_1 - A \uvec_2 \right\|_{V'} \le
\left( C(\overline{\alpha}) + C(\overline{\beta})
       \left( \|\uvec_1\|_{0,3,\Omega} + \|\uvec_2\|_{0,3,\Omega} \right) \right)
\| \uvec_1 - \uvec_2 \|_{0,3,\Omega} \: , $$
where $C(\overline{\alpha})$ and $C(\overline{\beta})$ are exactly 
the same constants as in (\ref{norm(Au)}).
Using inequality (\ref{ineq.2}) we obtain
\begin{equation}				\label{A:mon}
   \langle A\uvec - A\vvec , \uvec - \vvec \rangle_{V' \times V} \ge
   \frac{C(\underline{\beta})}{2} \|\uvec -\vvec\|_{0,3,\Omega}^3 \quad 
   \mbox{for} \quad \uvec, \vvec \in \left( \xLn{3}(\Omega) \right)^n \: ,
\end{equation}
where $C(\underline{\beta})$ depends only on $\underline{\beta}$ and $\Omega$,
too.
Therefore $A$ is strictly monotone for $V = \xWn{3}(\Div;\Omega)$, 
and uniformly monotone (and thus coercive) for $V = \left( \xLn{3}(\Omega) \right)^n$.
\end{proof}

\begin{rmrk}					\label{rem:hom.stat.pr}
Since $A$ is uniformly monotone for  $V = \left( \xLn{3}(\Omega) \right)^n$
we can conclude easily that $A$ is uniformly monotone (and thus coercive) for
$V =  \xWn{3}_0(\Div; \Omega) := \left\{ \vvec \in \xWn{3}(\Div; \Omega) \bigm|
				         \Div(\vvec) = 0 \right\}$.
Obviously, the solution $\mvec$ of the homogeneous problem (i.e., $f \equiv 0$)
satisfies $\mvec \in \xWn{3}_0(\Div; \Omega)$.
Therefore, in the homogeneous case, we can extend directly 
the proof of the unique solvability for the linear problem 
(\cf \cite[Prop.~I.1.1 and Thm.~I.1.1]{Brezzi/Fortin})
to the nonlinear problem (\ref{mix.form.stat.pr}).
Owing to the uniform monotonicity of $A$ we can use the theorem of Browder 
and Minty \cite[Thm.~26.A]{Zeidler} to show that there exists a unique solution
$\mvec \in \xWn{3}_0(\Div; \Omega)$ of (\ref{mix.form.stat.pr}).
The existence and uniqueness of a solution $S$ then follows exactly like
in the proof of \cite[Thm.~I.1.1]{Brezzi/Fortin}.
\end{rmrk}

%
\subsection{Regularization of the stationary problem}
For general $f$ the situation is not so simple as depicted 
in Remark \ref{rem:hom.stat.pr}.
We use regularization to show the existence
of a solution $(\mvec,S) \in V \times Q$ to (\ref{mix.form.stat.pr}),
where $V := \xWn{3}(\Div; \Omega)$ and $Q := \xLn{{3/2}}(\Omega)$.
For $\varepsilon > 0$ we define nonlinear forms $d_{\varepsilon}: V \times V \to \xR$
and $c_\varepsilon: Q \times Q \to \xR$ by
$$ d_{\varepsilon} (\uvec,\vvec) := 
   \varepsilon \int_\Omega |\Div(\uvec)| \Div(\uvec) \Div(\vvec) \, d\xvec
   ~\mbox{ for } \uvec, \vvec \in \xWn{3}(\Div;\Omega) \;\; , \quad
   c_\varepsilon(p,q) := 
   \varepsilon \int_\Omega \frac{p}{\sqrt{|p|}} q \, d\xvec
   ~\mbox{ for } p,q \in \xLn{{3/2}}(\Omega) \: . $$
Then the regularized problem is:
Find $(\mvec_\varepsilon,S_\varepsilon) \in V \times Q$ such that
\begin{equation}			\label{reg.stat.pr}
\begin{array}{rcl}
a(\mvec_\varepsilon,\vvec) + d_\varepsilon(\mvec_\varepsilon,\vvec) - b(\vvec,S_\varepsilon)
& = & g(\vvec) \quad \mbox{for all } \vvec \in V \; , \\	
c_\varepsilon(S_\varepsilon,q) + b(\mvec_\varepsilon,q) & = & f(q)
\quad \mbox{for all } q \in Q \; .
\end{array}
\end{equation}
Analogously to the definition of $A$ we define operators
$A_\varepsilon : V \to V'$ by
$ \langle A_\varepsilon \uvec , \vvec \rangle_{V' \times V} := 
  a(\uvec,\vvec) + d_{\varepsilon} (\uvec,\vvec)$ 
and $C_\varepsilon : Q \to Q'$ by
$ \langle C_\varepsilon p,q \rangle_{Q' \times Q} := c_\varepsilon(p,q)\,$.
It is evident that $A_\varepsilon \uvec$ is a linear functional on $V$
for all $\uvec \in V$, 
and $C_\varepsilon p$ is a linear functional on $Q$ for all $p \in Q$.
Again, continuity of $A_\varepsilon \uvec$ and $C_\varepsilon p$, resp.,
follow from the boundedness, which is obtained by means of
H\"older's inequality:
$$ \begin{array}{rcl}
   a(\uvec,\vvec) + d_{\varepsilon} (\uvec,\vvec) & \le &
   \left( C(\overline{\alpha}) \| \uvec \|_{0,3,\Omega} 
         + C(\overline{\beta}) \| \uvec\|_{0,3,\Omega}^2 \right)
   \| \vvec\|_{0,3,\Omega}
   + \varepsilon \|\Div(\uvec)\|_{0,3,\Omega}^2 \|\Div(\vvec)\|_{0,3,\Omega} \: , \\
   \left| c_\varepsilon(p,q) \right| & \le &
   \varepsilon \|p\|_{0,3/2,\Omega}^{1/2} \|q\|_{0,3/2,\Omega} \: . 
\end{array} $$

To show that (\ref{reg.stat.pr}) has a solution, we need
continuity, coercivity and monotonicity of $A_\varepsilon$ and $C_\varepsilon$.
For $A_\varepsilon$ these properties are consequences of (\ref{ineq.1}-\ref{ineq.2})
again, for $C_\varepsilon$ we need an additional lemma:
\begin{lmm}					\label{lemma:ineq.3&4}
The real-valued function $f: \xR \to \xR$, $x \mapsto |x|^{-1/2} x$
is strictly monotone and H\"older continuous of order $1/2$ on $\xR$
with H\"older constant $\sqrt{2}$, i.e.\ for all $x,y \in \xR$ it holds
\begin{equation}				\label{ineq.3}
\left| |x|^{-1/2} x - |y|^{-1/2} y \right| \le \sqrt{2} |x-y|^{1/2} \: .
\end{equation}
Furthermore
\begin{equation}				\label{ineq.4}
\frac{|x-y|^2}{\sqrt{|x|} + \sqrt{|y|}} \le
\left( \frac{x}{\sqrt{|x|}} - \frac{y}{\sqrt{|y|}} \right) (x-y) \: .
\end{equation}
\end{lmm}

\begin{prpstn} 					\label{prop:A_eps,C_eps}
For every $\varepsilon > 0$:
\begin{enumerate} \renewcommand{\labelenumi}{\alph{enumi})}
\item the operator $A_\varepsilon : V \to V'$ is continuous, coercive 
      and strictly monotone on $V$,
\item the operator $C_\varepsilon : Q \to Q'$ is continuous, coercive
      and strictly monotone on $Q$.
\end{enumerate}
\end{prpstn}
\begin{proof}
Ad a): ~ The continuity of $A_\varepsilon$ follows from Proposition~\ref{prop:A}
and the inequality
$$ \left| d_\varepsilon(\uvec_1,\vvec) - d_\varepsilon(\uvec_2,\vvec) \right| \le
   \varepsilon \left( \| \Div(\uvec_1) \|_{0,3,\Omega} 
	      + \| \Div(\uvec_2) \|_{0,3,\Omega} \right)
   \left\| \Div(\uvec_1) - \Div(\uvec_2) \right\|_{0,3,\Omega}
   \left\| \Div(\vvec) \right\|_{0,3,\Omega} $$
for all $\uvec_1,\uvec_2,\vvec \in V$, which is obtained by means 
of H\"older's inequality and of (\ref{ineq.1}).
Furthermore, an application of (\ref{ineq.2}) yields
$$ d_\varepsilon(\uvec,\uvec-\vvec) - d_\varepsilon(\vvec,\uvec-\vvec) \ge
   \frac{\varepsilon}{2} \left\| \Div(\uvec - \vvec) \right\|_{0,3,\Omega}^3 
   \quad \mbox{for all } \uvec, \vvec \in V \: . $$
Together with (\ref{A:mon}) this inequality implies 
the uniform monotonicity of $A_\varepsilon$:
$$ \langle A_\varepsilon \uvec - A_\varepsilon \vvec , \uvec - \vvec \rangle_{V \times V'}
   \ge \frac{1}{2} \min \left( C(\underline{\beta}) , \varepsilon \right) 
       \| \uvec-\vvec \|_V^3 \quad \mbox{for all } \uvec, \vvec \in V \: . $$
The coercivity and strict monotonicity of $A_\varepsilon$ are direct consequences
of the uniform monotonicity.

Ad b): ~ Using H\"older's inequality and (\ref{ineq.3}) we obtain
the continuity of $C_\varepsilon$:
$$ \left\| C_\varepsilon p - C_\varepsilon q \right\|_{Q'} \le
   \varepsilon \left\| \frac{p}{\sqrt{|p|}} - \frac{q}{\sqrt{|q|}} 
                \right\|_{0,3,\Omega} \le
   \sqrt{2} \, \varepsilon \left\| p - q \right\|_{0,3/2,\Omega}^{1/2}
   \quad \mbox{for all } p, q \in Q \: . $$
The coercivity follows from the equation
$$ \langle C_\varepsilon q , q \rangle_{Q' \times Q} = 
   \int_\Omega \varepsilon \frac{q^2}{\sqrt{|q|}} \, d\xvec =
   \varepsilon \|q\|_{0,3/2,\Omega}^{3/2} \: , $$
and strict monotonicity 
from the strict monotonicity of $x \mapsto |x|^{-1/2} x$:
$$ \langle C_\varepsilon p - C_\varepsilon q , p - q \rangle_{Q' \times Q} = 
   \varepsilon \int_\Omega 
   \left( \frac{p}{\sqrt{|p|}} - \frac{q}{\sqrt{|q|}} \right) (p-q) \, d\xvec
   > 0 \quad \mbox{for all } p,q \in Q \mbox{ with } p \neq q \: . $$
\end{proof}

Now we are in a position to prove:
\begin{prpstn}				\label{prop:reg.stat.pr}
For every $\varepsilon > 0$ there is a unique solution 
$(\mvec_\varepsilon,S_\varepsilon) \in V \times Q$ 
of the regularized problem \eqref{reg.stat.pr}.
\end{prpstn}
\begin{proof}
Adding the left hand sides of (\ref{reg.stat.pr}) we obtain
the following nonlinear form defined on $V \times Q$:
$$ \mathbf{a}_\varepsilon \! \left( (\uvec,p),(\vvec,q) \right) := a(\uvec,\vvec) 
   + d_\varepsilon(\uvec,\vvec) - b(\vvec,p) + c_\varepsilon(p,q) + b(\uvec,q) 
   \quad \mbox{for } (\uvec,p), (\vvec,q) \in V \times Q $$
and a nonlinear operator $\mathcal{A}_\varepsilon : (V \times Q) \to (V \times Q)'$
defined by
$ \langle \mathcal{A}_\varepsilon (\uvec,p) , (\vvec,q) \rangle_
   {(V \times Q)' \times (V \times Q)} = 
   \mathbf{a}_\varepsilon \! \left( (\uvec,p),(\vvec,q) \right) $.
To study the properties of $\mathcal{A}_\varepsilon$ we introduce 
the continuous linear operator $B: V \to Q'$ 
and its adjoint operator $B': Q \to V'$, defined by
$\langle B\vvec,q \rangle_{Q \times Q'} = b(\vvec,q) = 
 \langle B'q,\vvec \rangle_{V \times V'}$.
Then we can write 
$$ \langle \mathcal{A}_\varepsilon (\uvec,p) , (\vvec,q) \rangle_
   {(V \times Q)' \times (V \times Q)} =
   \langle A_\varepsilon \uvec , \vvec \rangle_{V' \times V}
   - \langle B'p , \vvec \rangle_{V' \times V}
   + \langle C_\varepsilon p , q \rangle_{Q' \times Q} + \langle B \uvec, q \rangle_{Q' \times Q} $$ 
for $(\uvec,p), (\vvec,q) \in V \times Q$.
Since $A_\varepsilon$, $C_\varepsilon$, $B$ and $B'$ are continuous,
$\mathcal{A}_\varepsilon$ is continuous, too.
Furthermore, $\mathcal{A}_\varepsilon$ is coercive and strictly monotone.
This is a straightforward consequence of the corresponding properties 
of $A_\varepsilon$ and $C_\varepsilon$, since the terms containing $B$ or $B'$ 
cancel each other.
For the strict monotonicity this reads
\begin{eqnarray*}
\lefteqn{\langle \mathcal{A}_\varepsilon (\uvec,p) , (\uvec-\vvec,p-q) 
             \rangle_{(V \times Q)' \times (V \times Q)} -
	 \langle \mathcal{A}_\varepsilon (\vvec,q) , (\uvec-\vvec,p-q) 
   	     \rangle_{(V \times Q)' \times (V \times Q)}} \hspace{0em} & & \\
& = & \langle A_\varepsilon \uvec - A_\varepsilon \vvec , \uvec-\vvec \rangle_{V' \times V}
    - \langle B' (p-q) , \uvec-\vvec \rangle_{V' \times V}
    + \langle C_\varepsilon p - C_\varepsilon q , p-q \rangle_{Q' \times Q}
    + \langle B (\uvec-\vvec), p-q \rangle_{Q' \times Q} \\
& = & \langle A_\varepsilon \uvec - A_\varepsilon \vvec , \uvec-\vvec \rangle_{V' \times V}
    - \langle B (\uvec-\vvec), p-q \rangle_{Q' \times Q}
    + \langle C_\varepsilon p - C_\varepsilon q , p-q \rangle_{Q' \times Q}
    + \langle B (\uvec-\vvec), p-q \rangle_{Q' \times Q} \\
& = & \langle A_\varepsilon \uvec - A_\varepsilon \vvec , \uvec-\vvec \rangle_{V' \times V}
    + \langle C_\varepsilon p - C_\varepsilon q , p-q \rangle_{Q' \times Q} > 0 \; ,
\end{eqnarray*}
if $\uvec \neq \vvec$ or $p \neq q$.
The proof of the coercivity of $\mathcal{A}_\varepsilon$ is even more simple.

Thus we can apply the theorem of Browder and Minty \cite[Thm.~26.A]{Zeidler} 
to show that for every $\mathbf{f} \in (V \times Q)'$ there exists 
a unique solution $(\mvec_\varepsilon , S_\varepsilon) \in V \times Q$ 
of the operator equation
$\mathcal{A}_\varepsilon (\mvec_\varepsilon , S_\varepsilon) = \mathbf{f} \,$.
In particular, we choose the linear form $\mathbf{f}$ defined by
$ \mathbf{f}(\vvec,q) := g(\vvec) + f(q) $,
which arises by adding the right hand sides of (\ref{reg.stat.pr}).
Therefore (\ref{reg.stat.pr}) has a unique solution.
\end{proof}

Next, we show that the solution $(\mvec_\varepsilon , S_\varepsilon)$ is bounded
independently of $\varepsilon$.

\begin{prpstn}				\label{prop:reg.stat.pr:est}
There exist constants $\mathcal{K}_{\mathbf{m}}, \mathcal{K}_S$, 
independent of $\varepsilon$, such that for sufficiently small $\varepsilon > 0$
the solution $(\mvec_\varepsilon, S_\varepsilon)$ of \eqref{reg.stat.pr}
satisfies the following estimates:
\begin{equation}			\label{reg.stat.pr:est}
\|\mvec_\varepsilon\|_V \le \mathcal{K}_{\mathbf{m}} \; , \quad
    \|S_\varepsilon\|_Q \le \mathcal{K}_S \; .
\end{equation}
\end{prpstn}

\begin{proof}
We begin with a bound for the norm of $\Div(\mvec_\varepsilon)$.
Using the second equation of (\ref{reg.stat.pr}) we obtain
$$ 
\|\Div(\mvec_\varepsilon)\|_{0,3,\Omega} = \|\Div(\mvec_\varepsilon)\|_{Q'} = 
\sup_{q \in Q} \frac{\left| b(\mvec_\varepsilon,q) \right|}{\|q\|_Q} = 
\sup_{q \in Q} \frac{\left| f(q) - c_\varepsilon(S_\varepsilon,q) \right|}{\|q\|_Q}
\le \| f \|_{0,3,\Omega} + \varepsilon \|S_\varepsilon\|_Q^{1/2} \; .
$$ 
The estimation of $\|\mvec_\varepsilon\|_{0,3,\Omega}$ is based on
the first equation in (\ref{reg.stat.pr}):
\begin{eqnarray*}
C(\underline{\beta}) \|\mvec_\varepsilon\|_{0,3,\Omega}^3 
& \le & \int_\Omega \beta |\mvec_\varepsilon| 
        \left( \mvec_\varepsilon \cdot \mvec_\varepsilon \right) \, d\xvec
\le a(\mvec_\varepsilon,\mvec_\varepsilon) 
   + d_\varepsilon (\mvec_\varepsilon,\mvec_\varepsilon)
= g(\mvec_\varepsilon) + b(\mvec_\varepsilon,S_\varepsilon) \\
& \le & \|g\|_{V'} \|\mvec_\varepsilon\|_{0,3,\Omega}
      + \left( \|g\|_{V'} + \|S_\varepsilon\|_Q \right)
        \|\Div(\mvec_\varepsilon)\|_{0,3,\Omega} \; .
\end{eqnarray*}
Together with the estimate for $\|\Div(\mvec_\varepsilon)\|_{0,3,\Omega}$ 
above this yields
\begin{equation}			\label{reg.stat.pr:m-L3-est}
\|\mvec_\varepsilon\|_{0,3,\Omega}^3 \le 
\frac{1}{C(\underline{\beta})}
\left( \|g\|_{V'} \|\mvec_\varepsilon\|_{0,3,\Omega}
      + \|g\|_{V'} \|f\|_{0,3,\Omega}
      + \varepsilon \|g\|_{V'} \|S_\varepsilon\|_Q^{1/2}
      + \|f\|_{0,3,\Omega} \|S_\varepsilon\|_Q 
      + \varepsilon \|S_\varepsilon\|_Q^{3/2} \right) \; .
\end{equation}
To bound $S_\varepsilon$ we employ the inf-sup condition (\ref{inf-sup}).
Together with the first equation in (\ref{reg.stat.pr}) and
the above estimate for $\|\Div(\mvec_\varepsilon)\|_{0,3,\Omega}$ we obtain 
\begin{eqnarray*}
\theta \|S_\varepsilon\|_Q
& \le & \sup_{\vvec \in V} \frac{b(\vvec,S_\varepsilon)}{\|\vvec\|_V}
 = \sup_{\vvec \in V} 
   \frac{a(\mvec_\varepsilon,\vvec) + d_\varepsilon(\mvec_\varepsilon,\vvec) - g(\vvec)}
        {\|\vvec\|_V} \\
& \le & \|A \mvec_\varepsilon\|_{V'} 
             + \varepsilon \|\Div(\mvec_\varepsilon)\|_{0,3,\Omega}^2 + \|g\|_{V'}
  \le \|A \mvec_\varepsilon\|_{V'} + \|g\|_{V'}
              + \varepsilon \left( \|f\|_{0,3,\Omega} 
		           + \varepsilon \|S_\varepsilon\|_Q^{1/2} \right)^2
\end{eqnarray*}
for some constant $\theta>0$. Thus for sufficiently small $\varepsilon$ 
($\varepsilon < \theta^{1/3}$ is enough) it holds
$$ \|S_\varepsilon\|_Q \! \le \! \frac{1}{\theta-\varepsilon^3}
   \Big( \|A \mvec_\varepsilon\|_{V'} + \|g\|_{V'} 
	+ \varepsilon \|f\|_{0,3,\Omega}^2
	+ 2 \varepsilon^2 \|f\|_{0,3,\Omega} \|S_\varepsilon\|_Q^{1/2} \Big) \, , $$
such that
$$
\|S_\varepsilon\|_Q^{1/2} \le 
\left( \frac{2 \varepsilon^2}{\theta-\varepsilon^3} \|f\|_{0,3,\Omega} 
      + \left( \frac{1}{\theta-\varepsilon^3}
              \Big( \|A \mvec_\varepsilon\|_{V'} + \|g\|_{V'} +
	         \varepsilon \|f\|_{0,3,\Omega}^2 \Big) \right)^{1/2} \right) \: .
$$
An application of the estimate
$ \|A \mvec_\varepsilon\|_{V'} \le 
   C(\overline{\alpha}) \|\mvec_\varepsilon\|_{0,3,\Omega} +
   C(\overline{\beta}) \|\mvec_\varepsilon\|_{0,3,\Omega}^2 $
finally yields
\begin{equation}				\label{reg.stat.pr:S-est}
\|S_\varepsilon\|_Q^{1/2} 
\le \kappa_0 + \kappa_1 \|\mvec_\varepsilon\|_{0,3,\Omega}^{1/2}
   + \kappa_2 \|\mvec_\varepsilon\|_{0,3,\Omega} \: ,
\end{equation}
where the coefficients $\kappa_i$ are bounded,
independently of $\mvec_\varepsilon$ and $S_\varepsilon$,
for sufficiently small $\varepsilon$, e.g.\ $\varepsilon < \theta^{1/3}/2$,
$$ \kappa_0 :=
 \frac{2 \varepsilon^2}{\theta-\varepsilon^3} \|f\|_{0,3,\Omega}
 + \left( \frac{1}{\theta-\varepsilon^3} 
   \left( \|g\|_{V'}
  	  + \varepsilon \|f\|_{0,3,\Omega}^2 \right) \right)^{1/2} \; , \quad
\kappa_1 := 
 \left( \frac{C(\overline{\alpha})}{\theta-\varepsilon^3} \right)^{1/2} \; , \quad
\kappa_2 := 
  \left( \frac{C(\overline{\beta})}{\theta-\varepsilon^3} \right)^{1/2} \; . $$
Inserting (\ref{reg.stat.pr:S-est}) into (\ref{reg.stat.pr:m-L3-est})
we obtain after some calculations the inequality
$$ \left( 1 - \varepsilon \frac{\kappa_2^2}{\underline{\beta} c_\ell^3} \right)
   \|\mvec_\varepsilon\|_{0,3,\Omega}^3 \le
   \sum_{i=0}^5 \lambda_i \|\mvec_\varepsilon\|_{0,3,\Omega}^{i/2} \; , $$
where the $\lambda_i$ ($i=0,\ldots,5$) are independent of
$\|\mvec_\varepsilon\|_{0,3,\Omega}$ and bounded for sufficiently small
$\varepsilon < \overline{\varepsilon} \le \theta^{1/3}/2$. 
Thus there exists a constant $\mathcal{K}_1 < \infty$, 
independent of $\varepsilon$, such that
$\|\mvec_\varepsilon\|_{0,3,\Omega} \le \mathcal{K}_1$
for $\varepsilon \le \overline{\varepsilon}$.
Inserting this estimate into (\ref{reg.stat.pr:S-est}) we obtain
the bound for $\|S_\varepsilon\|_Q$ and using the above estimate for 
$\|\Div(\mvec_\varepsilon)\|_{0,3,\Omega}$ finally yields 
the bound for $\|\mvec_\varepsilon\|_V$.
\end{proof}

%
\subsection{Solvability of the stationary problem (\ref{mix.form.stat.pr})}
\begin{thrm}					\label{thm:stat.pr}
The mixed formulation \eqref{mix.form.stat.pr} of the stationary problem
\eqref{stat.pr} possesses a unique solution
$(\mvec,S) \in \xWn{3}(\Div; \Omega) \times \xLn{{3/2}}(\Omega)$.
\end{thrm}
\begin{proof}
Analogously to the definition of $\mathbf{a}_\varepsilon$
we add the left hand sides of (\ref{mix.form.stat.pr}) and obtain
the nonlinear form $\mathbf{a}$ defined by
$\mathbf{a} \left( (\uvec,p),(\vvec,q) \right) := 
 a(\uvec,\vvec) - b(\vvec,p) + b(\uvec,q) $
and the nonlinear operator 
$\mathcal{A} : (V \times Q) \to (V \times Q)'$
defined by
$ \langle \mathcal{A} (\uvec,p) , (\vvec,q) \rangle_
   {(V \times Q)' \times (V \times Q)} = 
   \mathbf{a} \left( (\uvec,p),(\vvec,q) \right) $.
Setting $\varepsilon = 1/n$ let $(\mvec_n,S_n)$ be the unique solution of the
regularized problem (\ref{reg.stat.pr}).
Since $\left( (\mvec_n,S_n) \right)_{n \in \xN}$ is a bounded sequence
in $V \times Q$, there exists a weakly convergent subsequence, 
again denoted by $\left( (\mvec_n,S_n) \right)_{n \in \xN}$,
with (weak) limit $(\mvec,S) \in V \times Q$. 
As
\begin{eqnarray*}
\left\| \mathcal{A}(\mvec_n,S_n) 
                 - \mathbf{f} \right\|_{(V \times Q)'}
& = & \sup_{0 \neq (\vvec,q) \in V \times Q} 
\frac{ \left| \mathbf{a} \left( (\mvec_n,S_n) , (\vvec,q) \right)
             - \mathbf{f}(\vvec,q) \right|}
     {\|(\vvec,q)\|_{V \times Q}} \\
& = & \sup_{(\vvec,q)}
\frac{ \left| \mathbf{a}_{1/n} \left( (\mvec_n,S_n) , (\vvec,q) \right)
              - d_{1/n}(\mvec_n,\vvec) - c_{1/n}(S_n,q) 
              - \mathbf{f}(\vvec,q) \right|}
     {\|(\vvec,q)\|_{V \times Q}} \\
& \le & \frac{1}{n} \left( \|\Div(\mvec_n)\|_{0,3,\Omega}^2
			  + \| S_n \|_Q^{1/2} \right) 
\stackrel{n \to \infty}{\longrightarrow} 0
\end{eqnarray*}
the sequence $\left( \mathcal{A}(\mvec_n,S_n) \right)_{n \in \xN}$
converges strongly in $V'$ to $\mathbf{f}$ defined by
$\mathbf{f}(\vvec,q) := g(\vvec) + f(q)$.
Thus we can conclude that $\mathcal{A} (\mvec,S) = \mathbf{f}$ 
in $(V \times Q)'$ (see e.g.\ \cite[p.~474]{Zeidler}), i.e.,
$(\mvec,S)$ is a solution of (\ref{mix.form.stat.pr}).

To show uniqueness, we consider two solutions
$(\mvec_1,S_1)$ and $(\mvec_2,S_2)$ of (\ref{mix.form.stat.pr}).
Using the test functions $\vvec = \mvec_1 - \mvec_2$
and $q = S_1 - S_2$ we obtain
\begin{eqnarray*}
a(\mvec_1,\mvec_1-\mvec_2) - a(\mvec_2,\mvec_1-\mvec_2)
- b(\mvec_1-\mvec_2,S_1) + b(\mvec_1-\mvec_2,S_2) & = & 0 \\
b(\mvec_1,S_1-S_2) - b(\mvec_2,S_1-S_2) & = & 0 \; .
\end{eqnarray*}
Adding these equations yields
$$ 0 = a(\mvec_1,\mvec_1-\mvec_2) - a(\mvec_2,\mvec_1-\mvec_2)
     = \langle A \mvec_1 - A \mvec_2 , \mvec_1-\mvec_2 \rangle_
       {(V \times Q)' \times (V \times Q)} \: . $$
Since $A$ is strictly monotone it follows $\mvec_1-\mvec_2=0$. 

If $\mvec \in V$ is given, $S \in Q$ is defined as solution of 
the variational equation
$b(\vvec,S) = g(\vvec) - a(\mvec,\vvec)$ for all $\vvec \in V$.
Therefore the uniqueness of $S$ is a direct consequence of the injectivity
of the operator $B': Q \to V'$, \cf \cite[{\S}II, Rem.~1.6]{Brezzi/Fortin}.
\end{proof}

%
\section{The semi-discrete problem}
\label{sec:semid.pr}
%
We return to the transient problem governed by (\ref{Da-Fo.eq}) 
and (\ref{cont.eq}).
We discretize (\ref{cont.eq}) in time using the implicit Euler method.
This yields not only a method to solve the transient problem numerically,
but also an approach to prove its solvability, 
the technique of semi-discretization.

We define a partition $0 = t^0 < t^1 < \ldots < t^K = T$ of the segment $(0,T)$
into $K$ intervals of constant length $\Delta t = T/K$, i.e.,
$t_k = k \Delta t$ for $k=0,\ldots,K$.
In the following for $k=0,\ldots,K$ we use the denotations
$S^k := S(\cdot,t^k)$ and $\mvec^k := \mvec(\cdot,t^k)$ for the unknown 
solutions and, analogously defined, $\alpha^k$, $\beta^k$ and $\gamma^k$ 
for the coefficient functions, $S_b^k$ for the boundary conditions 
and $f^k$ for the source term.
The initial condition $S(\cdot,t^0) = S^0(\cdot) \in W_0^{1,3/2}(\Omega)$
is given.

Using the equation of state $\rho=\rho(S)$ defined in (\ref{def:rho}),
the discretization in time of the continuity equation (\ref{cont.eq})
with the implicit Euler method yields for each $k \in \{1,\ldots,K\}$
\begin{equation}				\label{semid.pr}
\begin{array}{rcll} \displaystyle
 \left( \alpha^k(\xvec) + \beta^k(\xvec) |\mvec^k(\xvec)| \right) \mvec^k(\xvec)
+ \nabla S^k(\xvec) & = & 0 \: , & \xvec \in \Omega \: , \\
\displaystyle \frac{\phi(\xvec)}{\Delta t}
  \left( \gamma^k(\xvec) \frac{S^k(\xvec)}{\sqrt{|S^k(\xvec)|}}
        - \gamma^{k-1}(\xvec)
          \frac{S^{k-1}(\xvec)}{\sqrt{|S^{k-1}(\xvec)|}} \right)
+ \Div \left( \mvec^k(\xvec) \right) & = & f^k(\xvec) \: , \quad 
& \xvec \in \Omega \: , \\
S(\xvec) & = & S_b^k(\xvec) \: , & \xvec \in \partial\Omega \: .
\end{array}
\end{equation}
Note that for each $k \in \{1,\ldots,K\}$ the function $S^{k-1}$ is known.

For each $k \in \{1,\ldots,K\}$ we require $f^k \in \xLn{3}(\Omega)$, 
$S_b^k \in \mathrm{W}^{1/3,3/2}(\partial\Omega)$,
$\phi, \alpha^k, \beta^k, \gamma^k \in \xLinfty(\Omega)$ and additionally 
$$ \left. \begin{array}{l}
0 < \underline{\phi} \le \phi(\xvec) \le \overline{\phi} < \infty \: , \\
0 < \underline{\alpha} \le \alpha^k(\xvec) \le \overline{\alpha} < \infty \: , \\
0 < \underline{\beta} \le \beta^k(\xvec) \le \overline{\beta} < \infty \: , \\
0 < \underline{\gamma} \le \gamma^k(\xvec) \le \overline{\gamma} < \infty
\end{array} \right\}
\mbox{ for almost every } \xvec \in \Omega \; . $$

%
\subsection{Mixed formulation of the semi-discrete problem}
We continue to use the spaces $V = \xWn{3}(\Div;\Omega)$ and $Q = \xLn{{3/2}}(\Omega)$.
Then the variational formulation reads:
Find $(\mvec^k,S^k) \in V \times Q$ such that
\begin{equation} 			\label{var.form.semid.pr}
\hspace{-0.15em}
\begin{array}{r@{~\:}c@{~\:}l@{\quad}l}
\displaystyle \int_\Omega \! \left( \alpha^k + \beta^k |\mvec^k| \right)
               (\mvec^k \cdot \vvec) \, d\xvec
- \int_\Omega \! \Div(\vvec) S^k \, d\xvec & = & \displaystyle
- \int_{\partial\Omega} \! S_b^k (\vvec \cdot \nvec) \, d\sigma 
& \mbox{for all } \vvec \in V \: , \\
\displaystyle \int_\Omega \frac{\phi\, \gamma^k}{\Delta t}
            \frac{S^k}{\sqrt{|S^k|}} \, q \, d\xvec 
+ \int_\Omega \! \Div(\mvec^k) \, q \, d\xvec & = & \displaystyle
\int_\Omega f^k q \, d\xvec
+ \int_\Omega \frac{\phi \, \gamma^{k-1}}{\Delta t}
            \frac{S^{k-1}}{\sqrt{|S^{k-1}|}} \, q \, d\xvec 
& \mbox{for all } q \in Q \: .
\end{array}
\hspace{-0.25em}
\end{equation}

We introduce additional nonlinear forms $c^k$ on $Q \times Q$  defined by
$$ 
c^k(p,q) := 
   \int_\Omega \frac{\phi \, \gamma^k}{\Delta t}
   \frac{p}{\sqrt{|p|}} \, q \, d\xvec
$$ 
and  nonlinear operators $C^k : Q \to Q'$ by 
$\langle C^k p,q \rangle_{Q' \times Q} := c^k(p,q)$.
Again, it is evident that $C^k p$ is a linear mapping on $Q$ 
for all $p \in Q$.
The continuity of $C^k p$ is equivalent to its boundedness,
which in turn is a consequence of H\"older's inequality 
and the boundedness of $\phi$ and $\gamma$:
$$ \left| c^k(p,q) \right| \le 
   \frac{\overline{\phi} \, \overline{\gamma}}{\Delta t} \,
   \|p\|_{0,3/2,\Omega}^{1/2} \|q\|_{0,3/2,\Omega} \: . $$
In the same manner as in the proof of Proposition \ref{prop:A_eps,C_eps} b)
we obtain the continuity, coercivity and monotonicity of $C^k$.

\begin{prpstn}					\label{prop:C^k}
The operators $C^k : Q \to Q'$ are continuous, coercive and
strictly monotone on $Q$.
\end{prpstn}

We use $a$, $b$ and $g$ as defined in Section \ref{sec:stat.pr},
where $a=a^k$ and $g=g^k$ depend on $k$, 
because $\alpha$, $\beta$ and the boundary condition $S_b$ may change in time.
Then we can write the mixed formulation (\ref{var.form.semid.pr}) of (\ref{semid.pr})
in the following way:
Find $\left( \mvec^k,S^k \right) \in V \times Q$, such that
\begin{equation}				\label{mix.form.semid.pr}
\begin{array}{rcll}
a^k(\mvec^k,\vvec) - b(\vvec,S^k) & = & g^k(\vvec) & \mbox{for all } \vvec \in V \: , \\
c^k(S^k,q) + b(\mvec^k,q) & = & \tilde{f}^k(q) & \mbox{for all } q \in Q \; .
\end{array}
\end{equation}
Here $\tilde{f}^k \in Q'$ 
for $k=1,\ldots,K$ is defined by
$$ \tilde{f}^k(q) :=
   \int_\Omega \left( f^k + \frac{\phi \, \gamma^{k-1}}{\Delta t}
   \frac{S^{k-1}}{\sqrt{|S^{k-1}|}} \right) q \, d\xvec \; . $$
For the remainder of this section we restrict our considerations
to a fixed time step $k$. Thus we can omit the superscript $k$.

%
\subsection{Regularization of the semi-discrete problem}
We use the technique of regularization again.
Thus we consider, instead of (\ref{mix.form.semid.pr}), the following
regularized problem for $\varepsilon > 0$:
Find $(\mvec_\varepsilon,S_\varepsilon) \in V \times Q$ such that
\begin{equation}					\label{reg.semid.pr}
\begin{array}{rcl@{\quad}l}
a(\mvec_\varepsilon,\vvec) + d_\varepsilon(\mvec_\varepsilon,\vvec) - b(\vvec,S_\varepsilon)
& = & g(\vvec) & \mbox{for all } \vvec \in V \; , \\
c(S_\varepsilon,q) + b(\mvec_\varepsilon,q) & = & \tilde{f}(q)
& \mbox{for all } q \in Q \; .
\end{array}
\end{equation}
Here $d_\varepsilon (\uvec,\vvec) := 
   \varepsilon \int_\Omega |\Div(\uvec)| \Div(\uvec) \Div(\vvec) \, d\xvec$ 
is defined as in Section \ref{sec:stat.pr}.

In the same manner as Proposition \ref{prop:reg.stat.pr} we obtain:
\begin{prpstn}				\label{prop:reg.semid.pr}
For every $\varepsilon > 0$ there exists a unique solution 
$(\mvec_\varepsilon,S_\varepsilon) \in V \times Q$
of the regularized semi-discrete problem \eqref{reg.semid.pr}.
\end{prpstn}

Next, we show that the solution $(\mvec_\varepsilon , S_\varepsilon)$ of (\ref{reg.semid.pr})
is bounded independently of $\varepsilon$, too.
Since we added in (\ref{reg.semid.pr}) only one regularizing term $d_\varepsilon$,
we can use different techniques for the estimation of $\mvec_\varepsilon$ and $S_\varepsilon$.
In particular we obtain estimates that hold for every $\varepsilon>0$: 
\begin{prpstn}				\label{prop:reg.semid.pr:est}
There exist constants $\mathcal{K}_{\mathbf{m}}, \mathcal{K}_S$, 
independent of $\varepsilon$, such that the solution $(\mvec_\varepsilon, S_\varepsilon)$ 
of \eqref{reg.semid.pr} satisfies the following estimates:
\begin{equation}			\label{reg.semid.pr:est}
\|\mvec_\varepsilon\|_V \le \mathcal{K}_{\mathbf{m}} \; , \quad
    \|S_\varepsilon\|_Q \le \mathcal{K}_S \; .
\end{equation}
\end{prpstn}
\begin{proof}
As in the proof of Proposition \ref{prop:reg.stat.pr:est} we begin with
an estimate for the norm of $\Div(\mvec_\varepsilon)$:
$$ 
\|\Div(\mvec_\varepsilon)\|_{0,3,\Omega} 
\le \|\tilde{f}\|_{Q'} + 
    \frac{\overline{\phi} \, \overline{\gamma}}
	 {\Delta t} \|S_\varepsilon\|_Q^{1/2} \: .
$$ 
The estimation of $\|\mvec_\varepsilon\|_{0,3,\Omega}$ uses the following inequality,
established in the proof of Proposition \ref{prop:reg.stat.pr:est},
$$ C(\underline{\beta}) \|\mvec_\varepsilon\|_{0,3,\Omega}^3 
   \le \|g\|_{V'} \|\mvec_\varepsilon\|_{0,3,\Omega}
      + \left( \|g\|_{V'} + \|S_\varepsilon\|_Q \right)
        \|\Div(\mvec_\varepsilon)\|_{0,3,\Omega} $$
to derive
$$ \|\mvec_\varepsilon\|_{0,3,\Omega} \le
   \left( \frac{1}{C(\underline{\beta})} \|g\|_{V'} \right)^{1/2} \!\!
   + \left( \frac{1}{C(\underline{\beta})} 
       \Big( \|g\|_{V'} + \|S_\varepsilon\|_Q \Big)
       \|\Div(\mvec_\varepsilon)\|_{0,3,\Omega} \right)^{1/3} \, . $$

Together with the estimate for $\|\Div(\mvec_\varepsilon)\|_{0,3,\Omega}$
we obtain the following bound for $\|\mvec_\varepsilon\|_V$:
$$ 
\|\mvec_\varepsilon\|_V \le \kappa_1 + \kappa_2 \|S_\varepsilon\|_Q^{1/2} \: ,
$$ 
where the constants $\kappa_1$ and $\kappa_2$ are independent of $\varepsilon$ 
and $\|S_\varepsilon\|_Q$. 
To derive an estimation for $\|S_\varepsilon\|_Q$,
we use in (\ref{reg.semid.pr}) the test functions 
$\vvec = \mvec_\varepsilon$ and $q = S_\varepsilon$ and add the resulting equations.
Since $b$ is a bilinear form, we obtain the inequality
$$ c(S_\varepsilon,S_\varepsilon)
 \le a(\mvec_\varepsilon,\mvec_\varepsilon) 
 + d_\varepsilon(\mvec_\varepsilon,\mvec_\varepsilon)
 + c(S_\varepsilon,S_\varepsilon) 
 = g(\mvec_\varepsilon) + \tilde{f}(S_\varepsilon) \: . $$
Using the coercivity of $C$ and the bound for $\|\mvec_\varepsilon\|_V$ derived above
we can therefore conclude
\begin{eqnarray*}
\frac{\underline{\phi} \, \underline{\gamma}}{\Delta t}
\|S_\varepsilon\|_Q^{3/2} & \le & c(S_\varepsilon,S_\varepsilon)
\le g(\mvec_\varepsilon) + \tilde{f}(S_\varepsilon)
\le \|g\|_{V'} \|\mvec_\varepsilon\|_V + \|\tilde{f}\|_{Q'} \|S_\varepsilon\|_Q \\
& \le & \|g\|_{V'} \left( \kappa_1 + \kappa_2 
                          \left\| S_\varepsilon \right\|_Q^{1/2} \right)
+ \|\tilde{f}\|_{Q'} \|S_\varepsilon\|_Q \: .
\end{eqnarray*}
This yields the existence of a bound $\mathcal{K}_S$ for $\|S_\varepsilon\|_Q$. 
\end{proof}

%
\subsection{Solvability of the semi-discrete problem (\ref{mix.form.semid.pr})}
Again, we consider the limit $\varepsilon \to 0$
and obtain in the same manner as in Section \ref{sec:stat.pr} the existence
of a solution of the semi-discrete problem (\ref{mix.form.semid.pr}).
The proof of the uniqueness differs from the proof of Theorem \ref{thm:stat.pr}.

\begin{thrm}					\label{thm:semid.pr}
The mixed formulation \eqref{mix.form.semid.pr} of the semi-discrete problem
\eqref{semid.pr} possesses a unique solution
$(\mvec,S) \in \xWn{3}(\Div; \Omega) \times \xLn{{3/2}}(\Omega)$.
\end{thrm}
\begin{proof}
Like in the proof of Theorem \ref{thm:stat.pr} we add both equations
in (\ref{mix.form.semid.pr}) and obtain the nonlinear form $a$, 
defined on $(V \times Q) \times (V \times Q)$, 
and the linear form $\mathbf{f} \in (V \times Q)'$, defined by
$$ \mathbf{a} \big( (\uvec,p) , (\vvec,q) \big) :=
   a(\uvec,\vvec) - b(\vvec,p) + c(p,q) + b(\uvec,q) \quad , \quad
   \mathbf{f}(\vvec,q) := g(\vvec) + \tilde{f}(q) \: . $$
Again, the operator $\mathcal{A} : V \times Q \to (V \times Q)'$ is defined by
$\langle \mathcal{A} (\uvec,p) , (\vvec,q) \rangle_{(V \times Q)' \times (V \times Q)}
 = \mathbf{a} \big( (\uvec,p) , (\vvec,q) \big)$.
Choosing $\varepsilon = 1/n$ for $n \in \xN$ we obtain a sequence of unique solutions 
$(\mvec_n,S_n)$ of the regularized problems (\ref{reg.semid.pr}).
Owing to Proposition \ref{prop:reg.semid.pr:est} the sequence
$\left( (\mvec_n,S_n) \right)_{n \in \xN}$ is bounded in $V \times Q$. 
Thus there is a weakly convergent subsequence, again denoted by
$\left( (\mvec_n,S_n) \right)_{n \in \xN}$,
which converges to $(\mvec,S) \in V \times Q$.
In the same manner as in the proof of Theorem \ref{thm:stat.pr} we obtain the identity
$\mathcal{A} (\mvec,S) = \mathbf{f}$ in $(V \times Q)'$, i.e., $(\mvec,S)$ 
is a solution of the semi-discrete mixed formulation (\ref{mix.form.semid.pr}).

Uniqueness of the solution $(\mvec,S) \in V \times Q$ follows
from the strict monotonicity of $\mathcal{A}$, which, in turn,
is a consequence of the strict monotonicity of $A$ and $C$.
\end{proof}

%
\section{The transient problem}
\label{sec:trans.pr}
%
Finally, we address the continuous transient problem.
We restrict our considerations here to the case of homogeneous Dirichlet
boundary conditions.
Due to the lack of regularity of the solution $\mvec$, it is not possible
to handle more general boundary conditions as in the former sections.
Thus we consider the following problem:
\begin{equation}				\label{trans.pr}
\begin{array}{rcll} \displaystyle
\left( \alpha(\xvec,t) + \beta(\xvec,t) |\mvec(\xvec,t)| \right) 
\mvec(\xvec,t)
+ \nabla S(\xvec,t) & = & 0 \: , & (\xvec,t) \in \Omega \times [0,T] \: , \\[1ex]
\displaystyle \phi(\xvec) \frac{\partial \rho(S(\xvec,t),\xvec,t)}{\partial t}
+ \Div \left( \mvec(\xvec,t) \right) & = & f(\xvec,t) \: , \quad 
& (\xvec,t) \in \Omega \times [0,T] \: , \\[1ex]
S(\xvec,t) & = & 0 \: , & (\xvec,t) \in \partial\Omega \times [0,T] \: , \\[1ex]
S(\xvec,0) & = & S^0(\xvec) \: , & \xvec \in \Omega \: .
\end{array}
\end{equation}
Again, we require that $S^0 \in W_0^{1,3/2}(\Omega)$,
$\phi \in \xLinfty(\Omega)$ with lower and upper bound
$0 < \underline{\phi} \le \phi(\xvec) \le \overline{\phi} < \infty$
for almost every $\xvec \in \Omega$.
For every $t \in [0,T]$ the time-varying coefficient functions have to satisfy
the following assumptions:
$f(\cdot,t) \in \xLn{3}(\Omega)$ and
$\alpha(\cdot,t), \beta(\cdot,t), \gamma(\cdot,t) \in \xLinfty(\Omega)$ 
with upper and lower bounds 
$$ \left. \begin{array}{l}
0 < \underline{\alpha} \le \alpha(\xvec,t) \le \overline{\alpha}<\infty \: , \\
0 < \underline{\beta} \le \beta(\xvec,t) \le \overline{\beta} < \infty \: , \\
0 < \underline{\gamma} \le \gamma(\xvec,t) \le \overline{\gamma} < \infty
\end{array} \right\}
\mbox{ for almost every } \xvec \in \Omega 
\mbox{ and every } t \in [0,T] \; . $$
Furthermore, we require these coefficient functions to be Lipschitz continuous 
in time, i.e., there exist constants $L(\alpha)$, $L(\beta)$, $L(\gamma)$
and $L(f)$ such that for every $0 \le t_1 \le t_2 \le T$:
$$ \begin{array}{r@{~}c@{~}l}
  \left\| \alpha(t_1) - \alpha(t_2) \right\|_{0,\infty,\Omega} 
  & \le & L(\alpha) \, |t_1 - t_2| \: , \\
  \left\| \beta(t_1) - \beta(t_2) \right\|_{0,\infty,\Omega} 
  & \le & L(\beta) \, |t_1 - t_2| \: , \\
  \left\| \gamma(t_1) - \gamma(t_2) \right\|_{0,\infty,\Omega} 
  & \le & L(\gamma) \, |t_1 - t_2|
  \end{array}
  \quad \mbox{and} \quad
  \left\| f(t_1) - f(t_2) \right\|_{0,3,\Omega} \le L(f) \, |t_1 - t_2| \; .
$$

%
\subsection{A priori estimates for the solutions of the semi-discrete problems}
As mentioned above we use the technique of semi-discretization in time 
to show the existence of solutions of the transient problem (\ref{trans.pr}).
One important step has been done in Section \ref{sec:semid.pr}:
The existence and uniqueness of the solutions to the semi-discrete problems 
has been established.
In the next step, we have to consider the limit $\Delta t \to 0$ 
(or $K \to \infty$).
Similar to the regularization technique employed in the last two sections,
we therefore have to provide a priori estimates for the solutions 
of the semi-discrete problems, which are independent of $\Delta t$.
The bounds $\mathcal{K}_{\mathbf{m}}$ and $\mathcal{K}_S$ of
Proposition \ref{prop:reg.semid.pr:est} do not fulfill this requirement.
Thus we investigate the semi-discrete problem (\ref{var.form.semid.pr})
for homogeneous Dirichlet boundary conditions.
In a slightly different notation this problem reads:
\begin{equation}				\label{mix.form.hom.Dir.pr}
\begin{array}{rcll}
a^k(\mvec^k,\vvec) - b(\vvec,S^k) & = & 0 & \mbox{for all } \vvec \in V \: , \\
\displaystyle \int_\Omega \frac{\phi}{\Delta t} 
	\left( \rho^k(S^k)-\rho^{k-1}(S^{k-1}) \right) q \, d\xvec
+ b(\mvec^k,q) & = & f^k(q) & \mbox{for all } q \in Q \; ,
\end{array}
\end{equation}
where $\rho^k(S^k) := \gamma^k S^k / \sqrt{|S^k|}\:$.

\begin{lmm}				\label{lemma:S-est:hom.Dir.pr}
For sufficiently small $\Delta t > 0$ there exists a constant $C_S$,
independent of $\Delta t$ (and of $K$), such that
\begin{equation}			\label{S-est:hom.Dir.pr}
\| S^k \|_{0,3/2,\Omega} \le C_S \quad \mbox{for all }
0 \le k \le K \: .
\end{equation}
\end{lmm}

\begin{proof}
Choosing $\vvec=\mvec^k$ and $q=S^k$ in (\ref{mix.form.hom.Dir.pr}) 
and adding the resulting equations yields
\begin{equation}			\label{proof:S-est:hom.Dir.pr}
a^k(\mvec^k,\mvec^k) 
+ \int_\Omega \frac{\phi}{\Delta t}
   \left( \rho^k(S^k)-\rho^{k-1}(S^{k-1}) \right) S^k \, d\xvec = f^k(S^k) \: .
\end{equation}
Since $a^k(\mvec^k,\mvec^k) \ge 0$ this implies
$$ \int_\Omega \phi
       \left( \rho^k(S^k)-\rho^{k-1}(S^{k-1}) \right) S^k \, d\xvec \le
   \Delta t f^k(S^k) = \Delta t \int_\Omega f^k S^k \, d\xvec \: . $$
Estimating the right hand side using Young's inequality, we obtain:
$$ \int_\Omega f^k S^k \, d\xvec \le \int_\Omega |f^k| |S^k| \, d\xvec
   \le \int_\Omega \frac{1}{3} |f^k|^3 + \frac{2}{3} |S^k|^{3/2} \, d\xvec
   = \frac{1}{3} \|f^k\|_{0,3,\Omega}^3 
    + \frac{2}{3} \|S^k\|_{0,3/2,\Omega}^{3/2} \: . $$
In a similar manner we can treat the left hand side. Since
$$
\left| \int_\Omega \phi \, \rho^{k-1}(S^{k-1}) S^k \, d\xvec \right| \le 
\frac{1}{3} \int_\Omega \phi \, \gamma^{k-1} |S^{k-1}|^{3/2} \, d\xvec
+ \frac{2}{3} \int_\Omega \phi \, \gamma^{k-1} |S^k|^{3/2} \, d\xvec \: ,
$$
it follows that
$$
\int_\Omega \! \phi
         \left( \rho^k(S^k)-\rho^{k-1}(S^{k-1}) \right) S^k \, d\xvec \ge
\frac{1}{3} \! \int_\Omega \! \phi \, \gamma^k |S^k|^{3/2} \, d\xvec
- \frac{1}{3} \! \int_\Omega \! \phi \, \gamma^{k-1} |S^{k-1}|^{3/2} \, d\xvec
- \frac{2}{3} \! \int_\Omega \! \phi \left( \gamma^{k-1} - \gamma^k \right) 
			    |S^k|^{3/2} \, d\xvec \: .
$$
Due to the assumptions on $\gamma$ the integrand in the last term
can be bounded by
$$ \left| \phi \, \gamma^{k-1} - \phi \, \gamma^k \right| =
   \phi \left| \gamma^{k-1} - \gamma^k \right| \le
   \phi \, L(\gamma) \Delta t \le 
   \phi \, \gamma^k \, \frac{L(\gamma)}{\underline{\gamma}} \Delta t \; . $$
Merging all the above estimates together this results in
$$ \int_\Omega \phi \, \gamma^k |S^k|^{3/2} \, d\xvec
         - \int_\Omega \phi \, \gamma^{k-1} |S^{k-1}|^{3/2} \, d\xvec
   \le \Delta t \|f^k\|_{0,3,\Omega}^3 
 + 2 \left( \frac{1}{\underline{\phi} \, \underline{\gamma}} 
            + \frac{L(\gamma)}{\underline{\gamma}} \right) \Delta t 
   \int_\Omega \phi \, \gamma^k |S^k|^{3/2} \, d\xvec \: . $$
If $\Delta t$ is sufficient small such that
$C \Delta t := 2 \left( \frac{1}{\underline{\phi} \, \underline{\gamma}} 
                + \frac{L(\gamma)}{\underline{\gamma}} \right) \Delta t < 1$,
we can conclude that
$$ \int_\Omega \phi \, \gamma^k |S^k|^{3/2} \, d\xvec
   \le \frac{1}{1 - C \Delta t}
     \left( \int_\Omega \phi \, \gamma^{k-1} |S^{k-1}|^{3/2} \, d\xvec
           + \Delta t \|f^k\|_{0,3,\Omega}^3 \right) $$
for $k=1,\ldots,K$. 
As $1/(1 - C \Delta t) > 1$, we obtain by induction for all $k=0,\ldots,K$
\begin{eqnarray*}
\int_\Omega \phi \, \gamma^k |S^k|^{3/2} \, d\xvec
& \le & (1 - C \Delta t)^{-k}
        \left( \int_\Omega \phi \, \gamma^0 |S^0|^{3/2} \, d\xvec
              + \sum_{i=1}^k \Delta t \|f^i\|_{0,3,\Omega}^3 \right) \\
& \le & (1 - C \Delta t)^{-K}
        \left( \int_\Omega \phi \, \gamma^0 |S^0|^{3/2} \, d\xvec
              + T C_f^3 \right) \: ,
\end{eqnarray*}
where $C_f$ is an upper bound for $\|f\|_{0,3,\Omega}$.
Note that for $K \to \infty$ (i.e., $\Delta t = T/K \to 0$) 
the expression $(1 - C \Delta t)^{-K} = (1-CT/K)^{-K}$ 
tends to $e^{CT}$.
In particular, this expression remains bounded.
\end{proof}

\begin{lmm}				\label{lemma:m-est:hom.Dir.pr}
For sufficiently small $\Delta t$ there exists a constant $C_\mathbf{m}$,
independent of $\Delta t$ (and of $K$), such that
\begin{equation}			\label{m-est:hom.Dir.pr}
\| \mvec^k \|_{0,3,\Omega} \le C_\mathbf{m} \quad \mbox{for all }
0 \le k \le K \: .
\end{equation}
\end{lmm}
\begin{proof}
Choosing $q=S^k-S^{k-1}$ we obtain from the second equation 
in (\ref{mix.form.hom.Dir.pr})
$$ \int_\Omega \frac{\phi}{\Delta t} 
	\left( \rho^k(S^k)-\rho^{k-1}(S^{k-1}) \right) 
	\left( S^k-S^{k-1} \right) \, d\xvec
  + b(\mvec^k,S^k-S^{k-1}) = f^k (S^k-S^{k-1}) \: . $$
Since the first term is non-negative, this implies
$b(\mvec^k,S^k-S^{k-1}) \le f^k (S^k-S^{k-1})\,$.
Furthermore, we choose $\vvec=\mvec^k$ in the first equation 
of (\ref{mix.form.hom.Dir.pr}) belonging to time step
$k$ and $k-1$ and subtract the resulting equations.
Using the inequality above this yields
$$ 
a^k(\mvec^k,\mvec^k) - a^{k-1}(\mvec^{k-1},\mvec^k)
= b(\mvec^k,S^k) -  b(\mvec^k,S^{k-1}) \le f^k (S^k-S^{k-1}) \: .
$$ 
Again, we apply Young's inequality to show
\begin{eqnarray*}
\left| \int_\Omega \alpha^{k-1} 
	\left( \mvec^{k-1} \cdot \mvec^k \right) d\xvec \right| 
& \le & \frac{1}{2} \int_\Omega \alpha^k |\mvec^k|^2 \, d\xvec
+ \frac{1}{2} \int_\Omega \alpha^{k-1} |\mvec^{k-1}|^2 \, d\xvec
+ \frac{1}{2} \int_\Omega \left( \alpha^{k-1} - \alpha^k \right)
			  |\mvec^k|^2 \, d\xvec \: , \\
\left| \int_\Omega \beta^{k-1} |\mvec^{k-1}|
	        \left( \mvec^{k-1} \cdot \mvec^k \right) d\xvec \right| 
& \le & \frac{1}{3} \int_\Omega \beta^k |\mvec^k|^3 \, d\xvec
+ \frac{2}{3} \int_\Omega \beta^{k-1} |\mvec^{k-1}|^3 \, d\xvec
+ \frac{1}{3} \int_\Omega \left( \beta^{k-1} - \beta^k \right)
			  |\mvec^k|^3 \, d\xvec \: .
\end{eqnarray*}
Owing to the definition of $a$ we obtain
\begin{eqnarray*}
\frac{1}{2} \int_\Omega \alpha^k |\mvec^k|^2 \, d\xvec
- \frac{1}{2} \int_\Omega \alpha^{k-1} |\mvec^{k-1}|^2 \, d\xvec
- \frac{1}{2} \int_\Omega \left( \alpha^{k-1} - \alpha^k \right)
			  |\mvec^k|^2 \, d\xvec \quad & & \\
+ \: \frac{2}{3} \int_\Omega \beta^k |\mvec^k|^3 \, d\xvec
- \frac{2}{3} \int_\Omega \beta^{k-1} |\mvec^{k-1}|^3 \, d\xvec
- \frac{1}{3} \int_\Omega \left( \beta^{k-1} - \beta^k \right)
			  |\mvec^k|^3 \, d\xvec & & \\ 
\le a^k(\mvec^k,\mvec^k) - a^{k-1}(\mvec^{k-1},\mvec^k)
& \le & f^k (S^k-S^{k-1}) \: .
\end{eqnarray*}
Summing this relation for $i=1,\ldots,k$ yields
\begin{eqnarray*}
\frac{1}{2} \int_\Omega \! \alpha^k |\mvec^k|^2 \, d\xvec
+ \frac{2}{3} \int_\Omega \! \beta^k |\mvec^k|^3 \, d\xvec
& \!\!\! \le \!\!\! & \frac{1}{2} \int_\Omega \! \alpha^0 |\mvec^0|^2 \, d\xvec
       + \frac{2}{3} \int_\Omega \! \beta^0 |\mvec^0|^3 \, d\xvec \\
& \!\!\! & + \sum_{i=1}^k \left[ 
	\int_\Omega \frac{1}{2} \left( \alpha^{i-1} \! - \alpha^i \right)
		    |\mvec^i|^2 
		   + \frac{1}{3} \left( \beta^{i-1} \! - \beta^i \right)
		    |\mvec^i|^3 \, d\xvec
        + f^i (S^i \! - S^{i-1}) \right] .
\end{eqnarray*}
This inequality holds for $k=0,\ldots,K$
Due to (\ref{S-est:hom.Dir.pr}) and the Lipschitz-continuity of $f$,
the last term in this sum in bounded.
Indeed, for $k=1,\ldots,K$ it holds:
\begin{eqnarray*}
\left| \sum_{i=1}^k f^i (S^i-S^{i-1}) \right|
& = & f^k S^k - f^1 S^0 + \sum_{i=1}^{k-1} (f^i-f^{i+1}) S^i \\[-2ex]
& \le & \|f^k\|_{0,3,\Omega} \|S^k\|_{0,3/2,\Omega}
      + \|f^1\|_{0,3,\Omega} \|S^0\|_{0,3/2,\Omega}
      + \sum_{i=1}^{k-1} \|f^i-f^{i+1}\|_{0,3,\Omega} \|S^i\|_{0,3/2,\Omega} \\[-2ex]
& \le & 2 C_f C_S + \sum_{i=1}^{k-1} L(f) \Delta t C_S 
  \le \left( 2 C_f + T L(f) \right) C_S \: .
\end{eqnarray*}
Analogously, the estimation of the other terms is based on the 
Lipschitz-continuity of $\alpha$ and $\beta$:
\begin{eqnarray*}
\int_\Omega \frac{1}{2} \left( \alpha^{i-1} - \alpha^i \right) 
                        |\mvec^i|^2 
	+ \frac{1}{3} \left( \beta^{i-1} - \beta^i \right) |\mvec^i|^3\, d\xvec
& \!\!\! \le \!\!\! & \frac{1}{2} \frac{L(\alpha)}{\underline{\alpha}} \Delta t
    \!\int_\Omega\! \alpha^i |\mvec^i|^2  d\xvec
   + \frac{1}{3} \frac{L(\beta)}{\underline{\beta}} \Delta t
    \!\int_\Omega\! \beta^i |\mvec^i|^3  d\xvec \\
& \!\!\! \le \!\!\! & C(\alpha,\beta)
	\Delta t \!\int_\Omega\! \alpha^i |\mvec^i|^2 + \beta^i |\mvec^i|^3  d\xvec
 = C(\alpha,\beta) \Delta t \, a^i(\mvec^i,\mvec^i) \: ,
\end{eqnarray*}
where $C(\alpha,\beta) := 
\max \left\{ \frac{1}{2} \frac{L(\alpha)}{\underline{\alpha}} ,
	     \frac{1}{3} \frac{L(\beta)}{\underline{\beta}} \right\}$.
Summing up, we obtain for $k=1,\ldots,K$
$$ \sum_{i=1}^k \int_\Omega \frac{1}{2} \left( \alpha^{i-1} - \alpha^i \right)
		    |\mvec^i|^2 
		   + \frac{1}{3} \left( \beta^{i-1} - \beta^i \right)
		    |\mvec^i|^3 \, d\xvec
 \le C(\alpha,\beta) \sum_{i=1}^k \Delta t \, a^i(\mvec^i,\mvec^i)
 \le C(\alpha,\beta) \sum_{i=1}^K \Delta t \, a^i(\mvec^i,\mvec^i) \: . $$
Using (\ref{proof:S-est:hom.Dir.pr}) finally yields
\begin{eqnarray*}
\lefteqn{\sum_{k=1}^K \Delta t \, a^k(\mvec^k,\mvec^k) \, = \,
\sum_{k=1}^K \left[ \Delta t f^k (S^k) -
             \int_\Omega \phi \left( \rho^k(S^k) - \rho^{k-1}(S^{k-1}) \right)
				S^k \, d\xvec \right]} \qquad & & \\
& \!\! \le \!\! & \sum_{k=1}^K \Delta t \|f^k\|_{0,3,\Omega} \|S^k\|_{0,3/2,\Omega}
 - \sum_{k=1}^K \left[ 
   \frac{1}{3} \int_\Omega \! \phi \gamma^k |S^k|^{3/2} \, d\xvec
 - \frac{1}{3} \int_\Omega \! \phi \gamma^{k-1} |S^{k-1}|^{3/2} \, d\xvec
                \right] \\
& & + \, \sum_{k=1}^K \frac{2}{3} 
	\int_\Omega \phi \left( \gamma^{k-1} - \gamma^k \right) 
		|S^k|^{3/2} \, d\xvec \\
& \!\! \le \!\! & \sum_{k=1}^K \Delta t C_f C_S
 + \frac{1}{3} \int_\Omega \phi \, \gamma^0 |S^0|^{3/2} \, d\xvec
 - \frac{1}{3} \int_\Omega \phi \, \gamma^K |S^K|^{3/2} \, d\xvec
 + \frac{2}{3} \overline{\phi} \sum_{k=1}^K 
    \left\| \gamma^{k-1} - \gamma^k \right\|_{0,\infty,\Omega}
    \|S^k\|_{0,3/2,\Omega}^{3/2} \\
& \!\! \le \!\! & T C_f C_S 
  + \frac{2}{3} \overline{\phi} \, \overline{\gamma} C_S^{3/2}
    + \frac{2}{3} \overline{\phi} \sum_{k=1}^K L(\gamma) \Delta t \, C_S^{3/2} 
 \le T C_f C_S + \frac{2}{3} \overline{\phi} 
	\left( \overline{\gamma} + T L(\gamma) \right) C_S^{3/2} \: .
\end{eqnarray*}
Summarizing all the relations above, we obtain the following inequality,
which holds for $k=0,\ldots,K$
\begin{eqnarray*}
\frac{1}{2} \int_\Omega \alpha^k |\mvec^k|^2 \, d\xvec
+ \frac{2}{3} \int_\Omega \beta^k |\mvec^k|^3 \, d\xvec
& \!\! \le \!\! & \frac{1}{2} \int_\Omega \alpha^0 |\mvec^0|^2 \, d\xvec
 + \frac{2}{3} \int_\Omega \beta^0 |\mvec^0|^3 \, d\xvec \\
& & + \left( 2 C_f + T L(f) + C(\alpha,\beta) T C_f \right) C_S
+ C(\alpha,\beta) \frac{2}{3} \overline{\phi} 
	\left( \overline{\gamma} + T L(\gamma) \right) C_S^{3/2} \, .
\end{eqnarray*}
Since $\beta^k > \underline{\beta} > 0$, this yields the assertion. 
\end{proof}

\begin{lmm}				\label{lemma:S'-est:hom.Dir.pr}
For sufficiently small $\Delta t$ there exists a constant $C_{S'}$,
independent of $\Delta t$ (and of $K$), such that
\begin{equation}			\label{S'-est:hom.Dir.pr}
\sum_{k=1}^K \Delta t 
    \int_\Omega \left| \frac{S^k - S^{k-1}}{\Delta t} \right|^{3/2} d \xvec 
\le C_{S'} \: .
\end{equation}
\end{lmm}
\begin{proof}
In a similar manner as in the proof of Lemma \ref{lemma:m-est:hom.Dir.pr}
we obtain
$$ \int_\Omega \frac{\phi}{\Delta t} 
	\left( \rho^k(S^k)-\rho^{k-1}(S^{k-1}) \right) 
	\left( S^k-S^{k-1} \right) \, d\xvec
   + a^k(\mvec^k,\mvec^k) - a^{k-1}(\mvec^{k-1},\mvec^k) 
   = \int_\Omega f^k (S^k-S^{k-1}) \, d\xvec \: . $$
Summing up for $k=1,\ldots,K$ yields
\begin{eqnarray*}
\lefteqn{\!\sum_{k=1}^K \frac{\underline{\phi} \, \underline{\gamma}}{\Delta t} 
    		\int_\Omega \! \left( \! \frac{S^k}{\sqrt{|S^k|}} 
				  - \frac{S^{k-1}}{\sqrt{|S^{k-1}|}} \right) 
			    (S^k\!-S^{k-1}) \, d\xvec 
\le \sum_{k=1}^K \left[ \int_\Omega \! f^k (S^k\!-S^{k-1}) \, d\xvec
		- a^k(\mvec^k,\mvec^k) + a^{k-1}(\mvec^{k-1},\mvec^k) \right]}
\hspace{16em} & & \\
& \le \!\! & \left| \sum_{k=1}^K \int_\Omega f^k (S^k-S^{k-1}) \, d\xvec \right|
        + \sum_{k=1}^K \left( a^{k-1}(\mvec^{k-1},\mvec^k) - a^k(\mvec^k,\mvec^k) \right) \: .
\end{eqnarray*}
As we have seen in the proof of Lemma \ref{lemma:m-est:hom.Dir.pr},
the first term on the right hand side is bounded by 
$\left( 2 C_f + T L(f) \right) C_S$.
For the second term the following estimate holds
\begin{eqnarray*}
\lefteqn{\sum_{k=1}^K \left( a^{k-1}(\mvec^{k-1},\mvec^k)
			     - a^k(\mvec^k,\mvec^k) \right)} \hspace{5em} && \\
& \le \!\! & \sum_{k=1}^K \left(
       \frac{1}{2} \int_\Omega \alpha^{k-1} |\mvec^{k-1}|^2 \, d\xvec 
     - \frac{1}{2} \int_\Omega \alpha^k |\mvec^k|^2 \, d\xvec
     + \frac{1}{2} \int_\Omega \left( \alpha^{k-1} - \alpha^k \right)
			|\mvec^k|^2 \, d\xvec \right. \\
& & \rule{1.9em}{0ex} \left. 
     + \, \frac{2}{3} \int_\Omega \beta^{k-1} |\mvec^{k-1}|^3 \, d\xvec 
     - \frac{2}{3} \int_\Omega \beta^k |\mvec^k|^3 \, d\xvec
     + \frac{1}{3} \int_\Omega \left( \beta^{k-1} - \beta^k \right)
			|\mvec^k|^3 \, d\xvec \right) \\
& = \!\! & \frac{1}{2} \int_\Omega \alpha^0 |\mvec^0|^2 \, d\xvec
     - \frac{1}{2} \int_\Omega \alpha^K |\mvec^K|^2 \, d\xvec
     + \frac{2}{3} \int_\Omega \beta^0 |\mvec^0|^3 \, d\xvec 
     - \frac{2}{3} \int_\Omega \beta^K |\mvec^K|^3 \, d\xvec \\
& & \rule{0em}{0ex} + \sum_{k=1}^K \left(
      \frac{1}{2} \int_\Omega \left( \alpha^{k-1} - \alpha^k \right)
			|\mvec^k|^2 \, d\xvec 
     + \frac{1}{3} \int_\Omega \left( \beta^{k-1} - \beta^k \right)
			|\mvec^k|^3 \, d\xvec \right) \\
& \le \!\! & \rule{0em}{0ex} 
     \frac{1}{2} C(\overline{\alpha}) \|\mvec^0\|_{0,3,\Omega}^2
     + \frac{2}{3} C(\overline{\beta}) \|\mvec^0\|_{0,3,\Omega}^3
     + C(\alpha,\beta) \left( \! T C_f C_S + \frac{2}{3} \overline{\phi} 
    \left( \overline{\gamma} + T L(\gamma) \right) C_S^{3/2} \right) .
\end{eqnarray*}
Using (\ref{ineq.4}) we therefore showed that there exists a constant $C > 0$,
independent of $\Delta t$, such that
$$ \sum_{k=1}^K \Delta t 
	\int_\Omega \frac{1}{\sqrt{|S^k|}+\sqrt{|S^{k-1}|}}
		    \left( \frac{S^k-S^{k-1}}{\Delta t} \right)^2 d\xvec
   \le \sum_{k=1}^K \frac{1}{\Delta t} 
	\int_\Omega \left( \frac{S^k}{\sqrt{|S^k|}} 
			  - \frac{S^{k-1}}{\sqrt{|S^{k-1}|}} \right) 
		    \left( S^k-S^{k-1} \right) \, d\xvec
   \le C \: . $$
Applying H\"older's inequality finally yields the assertion:
\begin{eqnarray*}
\lefteqn{\sum_{k=1}^K \Delta t 
    	\int_\Omega \left| \frac{S^k - S^{k-1}}{\Delta t} \right|^{3/2} d \xvec
        = \sum_{k=1}^K \Delta t
	  \int_\Omega \left( \sqrt{|S^k|}+\sqrt{|S^{k-1}|} \right)^{3/4}
		\left( \frac{1}{\sqrt{|S^k|}+\sqrt{|S^{k-1}|}}
		       \left( \frac{S^k-S^{k-1}}{\Delta t} \right)^2 \right)^{3/4}
	d\xvec} \hspace{4em} & & \\
& \le \!\! & \sum_{k=1}^K \left[
	\left( \Delta t \int_\Omega \left( \sqrt{|S^k|}+\sqrt{|S^{k-1}|} \right)^3
			d\xvec \right)^{1/4}
        \left( \Delta t \int_\Omega \frac{1}{\sqrt{|S^k|}+\sqrt{|S^{k-1}|}}
			    	    \left( \frac{S^k-S^{k-1}}{\Delta t} \right)^2
			d\xvec \right)^{3/4} \right] \\
& \le \!\! & \left( \sum_{k=1}^K \Delta t 
		\int_\Omega \left( \sqrt{|S^k|}+\sqrt{|S^{k-1}|} \right)^3 d\xvec
 	\right)^{1/4} 
	\left( \sum_{k=1}^K \Delta t 
		\int_\Omega \frac{1}{\sqrt{|S^k|}+\sqrt{|S^{k-1}|}}
			\left( \frac{S^k-S^{k-1}}{\Delta t} \right)^2 d\xvec
 	\right)^{3/4} \\
& \le \!\! & T^{1/4} 2^{3/4} C_S^{3/8} C^{3/4} =: C_{S'} \: .
\end{eqnarray*}
\end{proof}

Next, we show that the mixed formulation (\ref{mix.form.hom.Dir.pr}) is
equivalent to a variational formulation of the time-discretized
parabolic equation (\ref{parabol.eq}).
To this end, we recall the nonlinear mapping $F$ of (\ref{def:F}).
For fixed time $t=t_k$, we define the nonlinear mapping 
$F^k: \left( \xLn{{3/2}}(\Omega) \right)^n \to \left( \xLn{3}(\Omega) \right)^n$
and its inverse $G^k$ defined by 
$G^k(\vvec) = \left( \alpha^k + \beta^k |\vvec| \right) \vvec$.
Note that 
$\int_\Omega G^k(\uvec) \cdot \vvec \, d\xvec = a^k(\uvec,\vvec)$
for $\uvec, \vvec \in \left( \xLn{3}(\Omega) \right)^n$.

\begin{prpstn}					\label{prop:equiv}
\begin{enumerate} \renewcommand{\labelenumi}{(\alph{enumi})}
\item If $S^k \in \xWn{{1,3/2}}(\Omega)$ is a solution 
    of the variational formulation: 
    Find $S^k \in \xWn{{1,3/2}}_0(\Omega)$ such that
    \begin{equation}				\label{var.form:hom.Dir.pr}
    \int_\Omega \frac{\phi}{\Delta t} 
	\left( \rho^k(S^k) - \rho^{k-1}(S^{k-1}) \right) q \, d\xvec +
    \int_\Omega F^k \left( \nabla S^k \right) \cdot \nabla q \, d\xvec =
    f^k (q) \quad \mbox{for all } q \in \xWn{{1,3/2}}_0(\Omega)\: ,
\end{equation}
    then $\left( F^k(\nabla S^k),S^k \right)$ is a solution 
    of the mixed formulation \eqref{mix.form.hom.Dir.pr}.
    In particular, $F^k(\nabla S^k) \in \xWn{3}(\Div;\Omega)$.
\item If $\left( \mvec^k,S^k \right) \in \xWn{3}(\Div;\Omega) \times \xLn{{3/2}}(\Omega)$
    is a solution of the mixed formulation \eqref{mix.form.hom.Dir.pr},
    then $S^k$ is a solution of the variational formulation \eqref{var.form:hom.Dir.pr}.
    In particular, $S^k \in \xWn{{1,3/2}}_0(\Omega)$.
\end{enumerate}
\end{prpstn}
\begin{proof}
Ad a) ~ Let $S^k$ be a solution of (\ref{var.form:hom.Dir.pr}).
We define $\mvec^k := F(\nabla S^k)$.
Then Green's formula yields
$$ \int_\Omega G^k(\mvec^k) \cdot \vvec \, d\xvec
   = \int_\Omega G^k(F^k(\nabla S^k)) \cdot \vvec \, d\xvec
   = \int_\Omega \nabla S^k \cdot \vvec \, d\xvec \\
   = - \int_\Omega \Div(\vvec) S^k \, d\xvec
   \quad \mbox{for all } \vvec \in \xWn{{3}}(\Div;\Omega) \: . $$
This is the first equation in (\ref{mix.form.hom.Dir.pr}).
To derive the second equation in (\ref{mix.form.hom.Dir.pr}), we consider 
(\ref{var.form:hom.Dir.pr}) for $q \in \mathcal{D}(\Omega) \subset W_0^{1,3/2}(\Omega)$,
and apply Green's formula again:
$$ \int_\Omega (\rho^k(S^k) - \rho^{k-1}(S^{k-1}) - f^k) \, q \, d\xvec
  = - \int_\Omega F^k(\nabla S^k) \cdot \nabla q \, d\xvec
  = - \int_\Omega \mvec^k \cdot \nabla q \, d\xvec \: . $$ 
Thus the difference $\rho^k(S^k) - \rho^{k-1}(S^{k-1}) - f \in \xLn{3}(\Omega)$
is the generalized divergence of $\mvec^k$; consequently
$\mvec^k \in \xWn{{3}}(\Div;\Omega)$.
Because $\mathcal{D}(\Omega)$ is densely embedded into $\xLn{{3/2}}(\Omega)$,
the second equation in (\ref{mix.form.hom.Dir.pr}) follows.

Ad b) ~ Now, let $(\mvec^k,S^k)$ be a solution of (\ref{mix.form.hom.Dir.pr}).
Green's formula then implies
$$ \int_\Omega G^k(\mvec^k) \cdot \vvec \, d\xvec
  = - \int_\Omega \Div(\vvec) S^k \, d\xvec \quad 
  \mbox{for all } \vvec \in \left( \mathcal{D}(\Omega) \right)^n \: . $$
Thus in the sense of distributions it holds 
$\nabla S^k = G^k(\mvec^k) \in \left( \xLn{{3/2}}(\Omega) \right)^n$.
Consequently, $S^k \in \xWn{{1,3/2}}(\Omega)$ and $\mvec^k = F^k(\nabla S^k)$.
To prove that $S^k$ fulfills (\ref{var.form:hom.Dir.pr}), we consider
$q \in \xWn{{1,3/2}}_0(\Omega) \subset \xLn{{3/2}}(\Omega)$ 
in the first equation of (\ref{mix.form.hom.Dir.pr}).
Another application of Green's formula yields
\begin{eqnarray*}
f^k(q) & = &
\int_\Omega (\rho^k(S^k) - \rho^{k-1}(S^{k-1})) \, q \, d\xvec
- \int_\Omega \Div(F^k(\nabla S^k)) q \, d\xvec \\
& = & \int_\Omega (\rho^k(S^k) - \rho^{k-1}(S^{k-1})) \, q \, d\xvec
     + \int_\Omega F^k(\nabla S^k) \cdot \nabla q \, d\xvec \: .
\end{eqnarray*}
Finally, we consider again the first equation of (\ref{mix.form.hom.Dir.pr})
for $\vvec \in \left( \mathcal{D}(\bar{\Omega}) \right)^n$.
Applying Green's formula we obtain
$$ 0 = \int_\Omega \nabla S^k \cdot \vvec \, d\xvec
      + \int_\Omega \Div(\vvec) S^k \, d\xvec =
   \int_{\partial\Omega} \gamma_0 S^k \, (\vvec \cdot \nvec) \, d\sigma \: . $$
Consequently, $\gamma_0 S^k = 0$ in $\xWn{{1/3,3/2}}(\partial\Omega)$,
i.e.\ $S^k \in \xWn{{1,3/2}}_0(\Omega)$.
\end{proof}

Using this equivalence, we obtain a bound for $S^k$
in the norm of $\xWn{{1,3/2}}(\Omega)$.
\begin{lmm}				\label{lemma:further.est:hom.Dir.pr}
For sufficiently small $\Delta t$ there exist constants $\overline{C}_{S}$,
$C_{\rho'}$ and $C_{\overline{\mathbf{m}}}$, all independent of $\Delta t$ 
(and of $K$), such that
\begin{eqnarray}			\label{2.S-est:hom.Dir.pr}
\| S^k \|_{1,3/2,\Omega} & \le & \overline{C}_S 
\quad \mbox{for all } 0 \le k \le K \: , \\	\label{rho'-est:hom.Dir.pr}
\left\| \frac{\rho^k(S^k) - \rho^{k-1}(S^{k-1})}{\Delta t} \right\|_{-1,3,\Omega}
& \le & C_{\rho'}
\quad \mbox{for all } 1 \le k \le K \: , \\	\label{div(m)-est:hom.Dir.pr}
\left\| \Div (\mvec^k) \right\|_{-1,3,\Omega} & \le & C_{\overline{\mathbf{m}}}
\quad \mbox{for all } 1 \le k \le K \: .
\end{eqnarray}
\end{lmm}
\begin{proof}
Proposition \ref{prop:equiv} gives $\nabla S^k = - G^k(\mvec^k)$.
Therefore (\ref{m-est:hom.Dir.pr}) implies
$$ \left\| \nabla S^k \right\|_{0,3/2,\Omega} = 
   \left\| G^k(\mvec^k) \right\|_{0,3/2,\Omega} \le
   C(\overline{\alpha}) C_{\mathbf{m}} + 
   C(\overline{\beta}) C_{\mathbf{m}}^2 =: C_G \: . $$
Together with (\ref{S-est:hom.Dir.pr}) we obtain (\ref{2.S-est:hom.Dir.pr}).

Also, this equivalence yields (\ref{rho'-est:hom.Dir.pr}), because by means of 
(\ref{var.form:hom.Dir.pr}) we have for all $q \in \xWn{{1,3/2}}_0(\Omega)$
\begin{eqnarray*}
\left| \int_\Omega \frac{\phi}{\Delta t} 
		\left( \rho^k(S^k) - \rho^{k-1}(S^{k-1}) \right) q \, d\xvec \right|
& = & \left| f^k(q) -
\int_\Omega F^k \left( \nabla S^k \right) \cdot \nabla q \,d\xvec \right|
 = \left| f^k(q) + \int_\Omega \mvec^k \cdot \nabla q \,d\xvec \right| \\
& \le & \| f^k \|_{0,3,\Omega} \|q\|_{0,3/2,\Omega} +
	\| \mvec^k \|_{0,3,\Omega} \| \nabla q \|_{0,3/2,\Omega} \\
& \le & \left( \| f^k \|_{0,3,\Omega} + \| \mvec^k \|_{0,3,\Omega} \right)
	\|q\|_{1,3/2,\Omega} \: .
\end{eqnarray*}
Finally, we obtain (\ref{div(m)-est:hom.Dir.pr}),
since the first equation of (\ref{mix.form.hom.Dir.pr}) yields
$$ \left| \int_\Omega \Div(\mvec^k) q \, d\xvec \right|
   = \left| f^k(q) -
	     \int_\Omega \frac{\phi}{\Delta t} 
		\left( \rho^k(S^k) - \rho^{k-1}(S^{k-1}) \right) q \, d\xvec \right|
 \le \left( \| f^k \|_{0,3,\Omega} + \overline{\phi} C_{\rho'} \right)
        \|q\|_{1,3/2,\Omega} \: . $$
for all $q \in \xWn{{1,3/2}}_0(\Omega)$
\end{proof}

%
\subsection{Solvability of the continuous problem}
Due to the existence of unique solutions to the semi-discrete 
mixed formulation (\ref{mix.form.hom.Dir.pr}) we obtain for every
$K \in \xN$ a $K+1$-tuple of solutions 
$\left( (\mvec_{\Delta t}^k,S_{\Delta t}^k) \right)_{k=0,\ldots,K} \in
 \left( \xWn{3}(\Div;\Omega) \times \xLn{{3/2}}(\Omega) \right)^{K+1}$.
Recall that $\Delta t = T/K$.
We denote these $K+1$-tuples with
$\mvec_{\Delta t} := (\mvec_{\Delta t}^k)_{k=0,\ldots,K} 
 \in \left( \xWn{3}(\Div;\Omega) \right)^{K+1}$
and $S_{\Delta t} := (S_{\Delta t}^k)_{k=0,\ldots,K} 
     \in \left( \xLn{{3/2}}(\Omega) \right)^{K+1}$.
We define step functions, e.g.\ 
$\Pi_{\Delta t} S_{\Delta t} \in \xLinfty(0,T;\xWn{{1,3/2}}_0(\Omega))$, 
which are piecewise constant in time, by
$$ \left( \Pi_{\Delta t} S_{\Delta t} \right)(t) := \left\{
   \begin{array}{l@{\; , \quad\mathrm{if}~}l}
	S_{\Delta t}^0 & t=0 \: , \\[0.5ex]
	S_{\Delta t}^k & (k-1) \Delta t < t \le k \Delta t \: , ~ 
			 k=1,\ldots,K \: ,
   \end{array} \right. $$
and piecewise linear (in time) functions
$\Lambda_{\Delta t} S_{\Delta t} \in \mathrm{C}([0,T];\xWn{{1,3/2}}_0(\Omega))$ fulfilling
$$ \left( \Lambda_{\Delta t} S_{\Delta t} \right)(t^k) = S_{\Delta t}^k
   \quad \mbox{for } k=0,\ldots,K \: . $$
The time derivative of $\Lambda_{\Delta t} S_{\Delta t}$ is
a piecewise constant step function with values
$$ \Lambda_{\Delta t}' S_{\Delta t} (t) :=
   \frac{\partial}{\partial t} 
   \left( \Lambda_{\Delta t} S_{\Delta t} \right)(t) =
   \frac{S_{\Delta t}^k - S_{\Delta t}^{k-1}}{\Delta t} \; ,
   \quad \mbox{if } (k-1) \Delta t < t < k \Delta t \: ,
   ~ k=1,\ldots,K \: . $$
In addition, we use piecewise constant approximations $\gamma_{\Delta t}$ and
$f_{\Delta t}$ of the coefficient functions $\gamma$ and $f$,
and piecewise constant operators $\rho_{\Delta t}$, $F_{\Delta t}$ and $G_{\Delta t}$.
Owing to the lemmas above the following bounds hold for sufficiently small
time step sizes $\Delta t$:
\begin{eqnarray*}
\left\| \Pi_{\Delta t} S_{\Delta t} 
	\right\|_{\xLinfty(0,T;\xWn{{1,3/2}}_0(\Omega))}
& \le & \overline{C}_S \: , \\
\left\| \Lambda_{\Delta t}' S_{\Delta t} 
	\right\|_{\xLn{{3/2}}(0,T;\xLn{{3/2}}(\Omega))}
& \le & C_{S'} \: , \\
\left\| \Pi_{\Delta t} \rho_{\Delta t} (S_{\Delta t})
	\right\|_{\xLinfty \left( 0,T;\xLn{3}(\Omega) \right)}
& \le & \overline{\gamma} \, \sqrt{C_S} \: , \\
\left\| \Lambda_{\Delta t}' \rho_{\Delta t} (S_{\Delta t}) 
	\right\|_{\xLinfty \left(0,T;\xWn{{-1,3}}(\Omega) \right)}
& \le & C_{\rho'} \: , \\ 
\left\| \Pi_{\Delta t} \mvec_{\Delta t}
	\right\|_{\xLinfty \left(0,T;(\xLn{3}(\Omega))^n \right)}
& \le & C_{\mathbf{m}} \: , \\
\left\| \Pi_{\Delta t} \Div(\mvec_{\Delta t}) 
	\right\|_{\xLinfty \left(0,T;\xWn{{-1,3}}(\Omega) \right)}
& \le & C_{\overline{\mathbf{m}}} \: , \\
\left\| \Pi_{\Delta t} G_{\Delta t} (\mvec_{\Delta t}) 
	\right\|_{\xLinfty(0,T;(\xLn{{3/2}}(\Omega))^n)}
& \le & C_G \: , \\
\left\| \rho^K(S^K) \right\|_{0,3,\Omega}
& \le & \overline{\gamma} \, \sqrt{C_S} \: .
\end{eqnarray*}
The third (and the last) inequality follow from
$$ \int_\Omega \left| \rho^k(S^k) \right|^3 \, d\xvec =
   \int_\Omega \left| \gamma^k \frac{S^k}{\sqrt{|S^k|}} \right|^3
               \, d\xvec \le
   \overline{\gamma}^3 \int_\Omega \left| S^k \right|^{3/2} \, d\xvec \le
   \overline{\gamma}^3 C_S^{3/2} \: . $$

Thus there exist subsequences, again indexed by $\Delta t$,
that converge in the corresponding weak*-topology; in detail
\begin{equation}				\label{weak*-conv}
\begin{array}{r@{~}c@{~}l@{\quad\mathrm{in}\quad}l}
\Pi_{\Delta t} S_{\Delta t} & \stackrel{*}{\rightharpoonup}
 & S & \xLinfty \big( 0,T;\xWn{{1,3/2}}_0(\Omega) \big) \: , \\
\Lambda_{\Delta t}' S_{\Delta t} & \rightharpoonup
 & S' & \xLn{{3/2}} \left(0,T;\xLn{{3/2}}(\Omega) \right) \: , \\
\Pi_{\Delta t} \rho_{\Delta t} (S_{\Delta t}) & \stackrel{*}{\rightharpoonup}
 & R & \xLinfty \left(0,T;\xLn{3}(\Omega) \right) \: , \\
\Lambda_{\Delta t}' \rho_{\Delta t} (S_{\Delta t}) & \stackrel{*}{\rightharpoonup}
 & R' & \xLinfty \left(0,T;\xWn{{-1,3}}(\Omega) \right) \: , \\
\Pi_{\Delta t} \mvec_{\Delta t} & \stackrel{*}{\rightharpoonup}
 & \mvec & \xLinfty \left(0,T;(\xLn{3}(\Omega))^n \right) \: , \\
\Pi_{\Delta t} \Div(\mvec_{\Delta t}) & \stackrel{*}{\rightharpoonup}
 & \overline{\mvec} & \xLinfty \left(0,T;\xWn{{-1,3}}(\Omega) \right) \: , \\
\Pi_{\Delta t} G_{\Delta t} (\mvec_{\Delta t}) & \stackrel{*}{\rightharpoonup}
 & \gvec & \xLinfty \left(0,T;(\xLn{{3/2}}(\Omega))^n \right) \: \\
\rho^K(S^K) & \rightharpoonup & R_{T} & \xLn{3}(\Omega) \: .
\end{array}
\end{equation}

\begin{prpstn}					\label{prop:R=rho(S)}
The limits $S$ of $\Pi_{\Delta t} S_{\Delta t}$ and 
$R$ of $\Pi_{\Delta t} \rho_{\Delta t} (S_{\Delta t})$ 
from \eqref{weak*-conv} satisfy $\rho(S) = R$ almost everywhere in
$(0,T) \times \Omega$.
\end{prpstn}
\begin{proof}
As well as $\Pi_{\Delta t} S_{\Delta t}$, also 
$\Lambda_{\Delta t} S_{\Delta t}$ is bounded in 
$\xLinfty \big(0,T;\xWn{{1,3/2}}_0(\Omega) \big)$.
In particular, $\Lambda_{\Delta t} S_{\Delta t}$
and its partial derivatives
$(\partial/\partial x_i) \Lambda_{\Delta t} S_{\Delta t}$
are bounded in $\xLn{{3/2}} \! \left( 0,T;\xLn{{3/2}}(\Omega) \right)$.
Owing to (\ref{S'-est:hom.Dir.pr})
$(\partial/\partial t) \Lambda_{\Delta t} S_{\Delta t} 
 = \Lambda_{\Delta t}' S_{\Delta t}$ is bounded, too, 
such that $\Lambda_{\Delta t} S_{\Delta t}$ is bounded 
in $\xWn{{1,3/2}} \left( (0,T) \times \Omega \right)$.
The Rellich--Kondrachov-Theorem yields that
$\xWn{{1,3/2}}$ is embedded compactly in $\xLn{{3/2}}$. 
Thus there exists a subsequence (again denoted by $\Lambda_{\Delta t}$), which
converges (strongly) in $\xLn{{3/2}} \left( (0,T) \times \Omega \right)$ to $S$.
Choosing a further subsequence, we obtain that 
$\Lambda_{\Delta t} S_{\Delta t}$ converges almost everywhere 
in $(0,T) \times \Omega$ to $S$.
Applying the mapping $\rho_{\Delta t}$ yields
$\lim_{\Delta t \to 0} \Lambda_{\Delta t} \rho_{\Delta t} (S_{\Delta t})
 = \rho(S)$ a.e.\ in $(0,T) \times \Omega$.
Since $\Lambda_{\Delta t} \rho_{\Delta t} (S_{\Delta t})$
is bounded in $\xLn{3} \! \left( (0,T) \times \Omega \right)$, 
we can conclude that
$\Lambda_{\Delta t} \rho_{\Delta t} (S_{\Delta t})$ weakly converges to 
$\rho(S)$ in $\xLn{3} \!\left( (0,T) \times \Omega \right)$, 
i.e.\ for all $q \in \xLn{{3/2}} \! \left( (0,T) \times \Omega \right)$ it holds 
$$ \lim_{\Delta t \to 0}
   \int_0^{T} \int_\Omega \Lambda_{\Delta t} \rho_{\Delta t} (S_{\Delta t}) \, q \, d\xvec \, dt
   = \int_0^{T} \int_\Omega \rho(S) \, q \, d\xvec \, dt \: . $$
On the other hand
$$ \lim_{\Delta t \to 0}
   \int_0^{T} \int_\Omega \Lambda_{\Delta t} \rho_{\Delta t} (S_{\Delta t}) q \, d\xvec \, dt
  = \lim_{\Delta t \to 0}
     \int_0^{T} \int_\Omega \Pi_{\Delta t} \rho_{\Delta t} (S_{\Delta t}) q \, d\xvec \, dt
  = \int_0^{T} \int_\Omega R q \, d\xvec \, dt \: . $$
Since the limit of a convergent sequence is unique, the assertion follows.
\end{proof}

For the remainder of this section, we denote by $\langle \cdot , \cdot \rangle$
the dual pairing between $\xWn{{-1,3}}(\Omega)$ and $\xWn{{1,3/2}}_0(\Omega)$.
\begin{prpstn}					\label{prop:S'=dS/dt}
\begin{enumerate}	\renewcommand{\labelenumi}{\alph{enumi})}
\item The identity $S' = (\partial/\partial t)S$ holds in the sense of 
    distributions from $(0,T)$ to $\xLn{{3/2}}(\Omega)$, 
    i.e., for all $\varphi \in \mathcal{D}((0,T))$ it holds:
    $$ \int_0^{T} S'(t) \varphi(t) \, dt = -\int_0^{T} S(t) \varphi'(t) \, dt
       \quad \mbox{in } \xLn{{3/2}}(\Omega) \: . $$
\item The identity $R' = (\partial/\partial t)R$ holds in the sense of 
    distributions from $(0,T)$ to $\xWn{{-1,3}}(\Omega)$,
    i.e., for all $\varphi \in \mathcal{D}((0,T))$ it holds:
    $$ \int_0^{T} R'(t) \varphi(t) \, dt = -\int_0^{T} R(t) \varphi'(t) \, dt
       \quad \mbox{in } \xWn{{-1,3}}(\Omega) \: . $$
\item The identity $\overline{\mvec} = \Div(\mvec)$ in the sense of 
    distributions on $\Omega$ holds almost everywhere in $(0,T)$, 
    i.e., for all $\Psi \in \mathcal{D}(\Omega)$ it holds:
    $$ \langle \overline{\mvec} , \Psi \rangle = 
       - \int_\Omega \mvec \cdot \nabla \Psi \, d\xvec
       \quad \mbox{a.e.\ in } ~ (0,T) \: . $$
\item The identity $\gvec = - \nabla S$ holds in 
    $\xLinfty \left(0,T;(\xLn{{3/2}}(\Omega))^n \right)$,
    i.e., for all $\vvec \in \xLone \left( 0,T; (\xLn{3}(\Omega))^n \right)$ it holds:
    $$ \int_0^{T} \int_\Omega \gvec \cdot \vvec \, d\xvec \, dt = 
      - \int_0^{T} \int_\Omega \nabla S \cdot \vvec \, d\xvec \, dt \: . $$
\end{enumerate}
\end{prpstn}

\begin{proof}
Ad a) ~ From the second equation in (\ref{weak*-conv}) we can conclude that
for all $\varphi \in \mathcal{D}((0,T))$
$$ \lim_{\Delta t \to 0} 
	\int_0^{T} \Lambda_{\Delta t}' S_{\Delta t}(t) \varphi(t) \, dt =
   \int_0^{T} S'(t) \varphi(t) \, dt \quad \mbox{in } \xLn{{3/2}}(\Omega)\: . $$
On the other hand, partial integration yields
\begin{eqnarray*}
\lefteqn{\int_0^{T} \Lambda_{\Delta t}' S_{\Delta t}(t) \varphi(t) \, dt
         = \sum_{k=1}^K \int_{(k-1) \Delta t}^{k \Delta t} 
			\Lambda_{\Delta t}' S_{\Delta t}(t) \varphi(t) \, dt}
& & \\
& = & \sum_{k=1}^K \Big(
	\Lambda_{\Delta t} S_{\Delta t}(k\Delta t) \varphi(k\Delta t) -
	\Lambda_{\Delta t} S_{\Delta t}((k-1)\Delta t) \varphi((k-1)\Delta t)
		    \Big)
     - \sum_{k=1}^K \int_{(k-1) \Delta t}^{k \Delta t} 
			\Lambda_{\Delta t} S_{\Delta t}(t) \varphi'(t) \, dt \\
& = & - \int_0^{T} \Lambda_{\Delta t} S_{\Delta t}(t) \varphi'(t) \, dt \: ,
\end{eqnarray*}
such that
\begin{eqnarray*}
\int_0^{T} S'(t) \varphi(t) \, dt & = & 
\lim_{\Delta t \to 0}
    \int_0^{T} \Lambda_{\Delta t}' S_{\Delta t}(t) \varphi(t) \, dt
 = \lim_{\Delta t \to 0} 
    	- \int_0^{T} \Lambda_{\Delta t} S_{\Delta t}(t) \varphi'(t) \, dt \\
& = & \lim_{\Delta t \to 0} 
    	- \int_0^{T} \Pi_{\Delta t} S_{\Delta t}(t) \varphi'(t) \, dt =
	- \int_0^{T} S(t) \varphi'(t) \, dt \: .
\end{eqnarray*}
The identity in b) follows in a similar manner as the identity in a).

Ad c) ~ Let $\Psi \in \mathcal{D}(\Omega)$ and $\varphi \in \mathcal{D}((0,T))$
be arbitrarily chosen. Then
\begin{eqnarray*}
\lim_{\Delta t \to 0} \int_0^{T}
    \langle \Pi_{\Delta t} \, \Div (\mvec_{\Delta t}) , \Psi \rangle \varphi(t) \, dt
& = & \lim_{\Delta t \to 0} \int_0^{T}
    \langle \Pi_{\Delta t} \, \Div (\mvec_{\Delta t}) , \varphi(t) \Psi \rangle \, dt \\
= \int_0^{T} \langle \overline{\mvec}(t) , \varphi(t) \Psi \rangle \, dt
& = & \int_0^{T} \langle \overline{\mvec}(t) , \Psi \rangle \varphi(t) \, dt \: ,
\end{eqnarray*}
because $\overline{\mvec}$ is the limit of 
$\left( \Pi_{\Delta t} \Div(\mvec_{\Delta t}) \right)_{\Delta t}$
in $\xLinfty \left(0,T;\xWn{{-1,3}}(\Omega) \right)$.
On the other hand
\begin{eqnarray*}
\int_0^{T} \langle \overline{\mvec}(t) , \Psi \rangle \varphi(t) \, dt
& = & \lim_{\Delta t \to 0} \int_0^{T}
    \langle \Pi_{\Delta t} \Div (\mvec_{\Delta t}) , \Psi \rangle \varphi(t) \, dt
 = \lim_{\Delta t \to 0}
    - \int_0^{T} \int_\Omega \Pi_{\Delta t} \mvec_{\Delta t} \cdot \nabla \Psi 
			\, d\xvec \, \varphi(t) \, dt \\
& = & \lim_{\Delta t \to 0} 
    - \int_0^{T} \int_\Omega \Pi_{\Delta t} \mvec_{\Delta t} \cdot \nabla \Psi 
			\, \varphi(t) \, d\xvec \, dt
= - \int_0^{T} \int_\Omega \mvec \cdot \nabla \Psi \, \varphi(t) \, d\xvec \, dt \\
& = & - \int_0^{T} \int_\Omega \mvec \cdot \nabla \Psi \, d\xvec \,
		 \varphi(t) \, dt \: . 
\end{eqnarray*}
Since $\varphi$ is arbitrarily chosen, the assertion follows.

Ad d) ~ We have seen in the proof of Proposition \ref{prop:equiv} that
$$ \int_\Omega G^k(\mvec^k) \cdot \vvec \, d\xvec =
   - \int_\Omega \nabla S^k \cdot \vvec \, d\xvec 
   \quad \mbox{for all } \vvec \in (\xLn{3} (\Omega))^n \: . $$
Consequently, for $\vvec \in \xLone \left( 0,T; (\xLn{3}(\Omega))^n \right)$ it holds
\begin{eqnarray*}
\int_0^{T} \int_\Omega \gvec \cdot \vvec \, d\xvec \, dt
& = & \lim_{\Delta t \to 0} 
    \int_0^{T} \int_\Omega \Pi_{\Delta t} G_{\Delta t}(\mvec_{\Delta t})
			 \cdot \vvec \, d\xvec  \, dt \\
& = & \lim_{\Delta t \to 0} - \int_0^{T} \int_\Omega 
	\nabla \big( \Pi_{\Delta t} S_{\Delta t} \big) \cdot \vvec \, d\xvec \, dt
  = - \int_0^{T} \int_\Omega \nabla S \cdot \vvec \, d\xvec \, dt \: .
\end{eqnarray*}
\end{proof}

The first two statements of Proposition \ref{prop:S'=dS/dt} imply that
after possible modification on a set of measure zero in $[0,T]$ we have
that $S \in \mathrm{C} \left( [0,T]; \xLn{{3/2}}(\Omega) \right)$ 
and $R \in \mathrm{C} \left( [0,T]; \xWn{{-1,3}}(\Omega) \right)$.
Thus for every $t \in [0,T]$ the value $S(t) \in \xLn{{3/2}}(\Omega)$
is well defined.

\begin{prpstn}				\label{prop:distr.var.form}
The following identity holds in 
$\xLinfty \! \left(0,T;\xWn{{-1,3}}(\Omega) \right)$:
$$ 
\phi \, \frac{\partial \rho(S)}{\partial t} + \Div(\mvec) = f \: .
$$ 
Furthermore $\rho(S(0)) = \rho^0(S^0)$ and $\rho(S(T)) = R_{T}\,$.
\end{prpstn}
   
\begin{proof}
For $\varphi \in \mathcal{D}(\overline{(0,T)})$ we define a step function
$\varphi_{\Delta t}$ by
$$ \varphi_{\Delta t}(t) := \left\{
   \begin{array}{l@{\; , \quad\mathrm{if}\quad}l}
	\varphi((k-1) \Delta t) & (k-1) \Delta t \le t < k \Delta t \: , ~ 
			 k=1,\ldots,K \: , \\[0.5ex]
	\varphi(T) & t=T \: .
   \end{array} \right. $$
Using the test function $q=\Psi \in \mathcal{D}(\Omega)$ in 
(\ref{mix.form.hom.Dir.pr}), 
multiplying by $\Delta t \, \varphi((k-1) \Delta t)$ 
and summing up for $k=1,\ldots,K$, we obtain
\begin{equation}			\label{sums}
\begin{array}{rcl}
\displaystyle \sum_{k=1}^K \Delta t
    \int_\Omega \phi \frac{\rho^k(S^k) - \rho^{k-1}(S^{k-1})}{\Delta t} 
		\Psi \, d\xvec \, \varphi((k-1) \Delta t) & & \\
\displaystyle + \sum_{k=1}^K \Delta t
	\int_\Omega \Div(\mvec^k) \, \Psi \, d\xvec \, \varphi((k-1) \Delta t)
& = & \displaystyle \sum_{k=1}^K \Delta t 
	\int_\Omega f^k \Psi \, d\xvec \, \varphi((k-1) \Delta t)
\end{array}
\end{equation}
Employing the piecewise constant functions $\Pi_{\Delta t}$
and $\Lambda_{\Delta t}'$ this reads
\begin{eqnarray*}
\sum_{k=1}^K \int_{(k-1) \Delta t}^{k \Delta t}
    \int_\Omega \phi \Lambda_{\Delta t}' \rho_{\Delta t} (S_{\Delta t})
		\Psi \, d\xvec \, \varphi_{\Delta t} \, dt \quad & & \\
+ \sum_{k=1}^K \int_{(k-1) \Delta t}^{k \Delta t}
	\int_\Omega \Pi_{\Delta t} \Div(\mvec_{\Delta t}) \, \Psi \, d\xvec 
	\varphi_{\Delta t} \, dt
& = & \sum_{k=1}^K \int_{(k-1) \Delta t}^{k \Delta t}
	\int_\Omega \Pi_{\Delta t} f_{\Delta t} \Psi \, d\xvec \,
	\varphi_{\Delta t} \, dt \: ,
\end{eqnarray*}
and after joining the integrals over $t$
$$ \int_0^{T} \int_\Omega \phi \Lambda_{\Delta t}' \rho_{\Delta t} (S_{\Delta t})
		\, \Psi \varphi_{\Delta t} \, d\xvec \, dt \quad
  + \int_0^{T} \int_\Omega \Pi_{\Delta t} \Div(\mvec_{\Delta t}) \, 
		\Psi \varphi_{\Delta t} \, d\xvec \, dt
  = \int_0^{T}
	\int_\Omega \Pi_{\Delta t} f_{\Delta t} 
		\, \Psi \varphi_{\Delta t} \, d\xvec \, dt \: . $$
Since $\left( \Psi \varphi_{\Delta t} \right)_{\Delta t}$ strongly converges 
to $\Psi \varphi$ in $\xLone \! \left( 0,T; \xWn{{1,3/2}}_0 (\Omega) \right)$,
we can pass to the limit $\Delta t \to 0$:
\begin{equation}				\label{ints}
\int_0^{T} \left\langle \phi \frac{\partial \rho(S)}{\partial t} , \Psi \varphi \right\rangle dt
 + \int_0^{T} \left\langle (\mvec) , \Psi \varphi \right\rangle dt
 = \int_0^{T} \int_\Omega f \, \Psi \varphi \, d\xvec \, dt \: .
\end{equation}
But the set
$\left\{ \Psi(\xvec) \varphi(t) \bigm| 
	 \Psi \in \mathcal{D}(\Omega) , \: \varphi \in \mathcal{D}(\overline{(0,T)}) \right\}$
is a dense subset of $\xLone \! \left( 0,T; \xWn{{1,3/2}}_0 (\Omega) \right)$.
Therefore the first identity in Proposition (\ref{prop:distr.var.form})
is established.

To prove the remaining two identities we first conclude from (\ref{sums}) that
\begin{eqnarray*}
\sum_{k=1}^K \left(
    \int_\Omega \phi \rho^k(S^k) \Psi \, d\xvec \, \varphi((k-1) \Delta t)
  - \int_\Omega \phi \rho^{k-1}(S^{k-1}) \Psi \, d\xvec \, \varphi((k-1) \Delta t)
  \right) & & \\
+ \sum_{k=1}^K \Delta t
	\int_\Omega \Div(\mvec^k) \, \Psi \, d\xvec \, \varphi((k-1) \Delta t)
& \!\! = \!\! & \sum_{k=1}^K \Delta t 
	\int_\Omega f^k \Psi \, d\xvec \, \varphi((k-1) \Delta t) \: .
\end{eqnarray*}
Rearranging the terms in the first line yields
\begin{eqnarray*}
\lefteqn{\sum_{k=1}^K \left(
    \int_\Omega \phi \rho^k(S^k) \Psi \, d\xvec \, \varphi((k-1) \Delta t)
  - \int_\Omega \phi \rho^{k-1}(S^{k-1}) \Psi \, d\xvec \, \varphi((k-1) \Delta t)
  	 \right)} \hspace{14.5em} & & \\
& = & - \sum_{k=1}^K \int_\Omega \phi \rho^k(S^k) \Psi \, d\xvec 
		\Big( \varphi(k \Delta t) - \varphi((k-1) \Delta t) \Big) \\
& & + \int_\Omega \phi \rho^K(S^K) \Psi \, d\xvec \, \varphi(K \Delta t)
    - \int_\Omega \phi \rho^0(S^0) \Psi \, d\xvec \, \varphi(0) \: .
\end{eqnarray*}
Like above this leads to
\begin{eqnarray*}
\lefteqn{- \int_0^{T} \int_\Omega \phi \Pi_{\Delta t} 
		\rho_{\Delta t} (S_{\Delta t}) \, \Psi \, d\xvec \, 
     	 \frac{\varphi(k \Delta t) - \varphi((k-1) \Delta t)}{\Delta t}\,dt
         + \int_0^{T} \int_\Omega \Pi_{\Delta t} \Div(\mvec_{\Delta t}) \, 
		\Psi \, d\xvec \, \varphi_{\Delta t} \, dt} \hspace{8em} & & \\
& = & \int_0^{T} \int_\Omega \Pi_{\Delta t} f_{\Delta t} \, 
	\Psi \, d\xvec \, \varphi_{\Delta t} \, dt
+ \int_\Omega \phi \rho^0(S^0) \Psi \, d\xvec \, \varphi(0)
- \int_\Omega \phi \rho^K(S^K) \Psi \, d\xvec \, \varphi(T) \: .
\end{eqnarray*}
Passing to the limit $\Delta t \to 0$ we obtain
$$ - \int_0^{T} \!\! \int_\Omega \phi \rho(S) \, \Psi \, d\xvec \,
		\frac{\partial \varphi}{\partial t}\,dt
	 + \int_0^{T} \!\! \langle(\mvec) , \Psi \rangle \varphi \, dt
  = \int_0^{T} \!\! \int_\Omega f \, \Psi \, d\xvec \, \varphi \, dt
   + \int_\Omega \phi \rho^0(S^0) \Psi \, d\xvec \, \varphi(0)
   - \int_\Omega \phi R_{T} \Psi \, d\xvec \, \varphi(T) \: . $$
In the other hand, partial integration of (\ref{ints}) yields
$$ 
- \int_0^{T} \!\!\! \int_\Omega \!\! \phi \rho(S) \, \Psi \, d\xvec \,
		\frac{\partial \varphi}{\partial t}\,dt
	 + \int_0^{T} \!\! \langle \Div(\mvec) , \Psi \rangle \varphi \, dt
 = \int_0^{T} \!\!\! \int_\Omega \! f \, \Psi \, d\xvec \, \varphi \, dt
+ \int_\Omega \!\! \phi \rho(S(0)) \Psi \, d\xvec \, \varphi(0)
- \int_\Omega \!\! \phi \rho(S(T)) \Psi \, d\xvec \, \varphi(T) \: .
$$ 
Subtracting the last two equations we can conclude
$$ \left( \int_\Omega \phi \rho^0(S^0) \Psi \, d\xvec
      - \int_\Omega \phi \rho(S(0)) \Psi \, d\xvec \right) \varphi(0)
 - \left( \int_\Omega \phi R_{T} \Psi \, d\xvec
- \int_\Omega \phi \rho(S(T)) \Psi \, d\xvec \right) \varphi(T) = 0 \: . $$
Since $\varphi(0)$ and $\varphi(T)$ are arbitrary, this implies
$$
\int_\Omega \phi \rho^0(S^0) \Psi \, d\xvec =
\int_\Omega \phi \rho(S(0)) \Psi \, d\xvec 
\quad \mbox{and} \quad
\int_\Omega \phi R_{T} \Psi \, d\xvec =
\int_\Omega \phi \rho(S(T)) \Psi \, d\xvec
$$
and finally $\rho(S(0)) = \rho^0(S^0)$ and $\rho(S(T)) = R_{T}$.
\end{proof}

Only the identity $\gvec = G(\mvec)$ is missing yet.
To show this, we need an auxiliary result, a generalization of
Lemma~1.2 from \cite{Raviart:70}.

\begin{lmm}					\label{lemma:int(dR/dt,S)}
The limit $S$ of $\left( \Pi_{\Delta t} S_{\Delta t} \right)_{\Delta t}$ 
satisfies:
$$ \int_0^{T} 
  \left\langle \phi \frac{\partial \rho(S)}{\partial t} , S \right\rangle \, dt
= \frac{1}{3} \left( \int_\Omega \phi \, |\rho(S(T))| \, |S(T)| \, d\xvec
	- \int_\Omega \phi \, |\rho(S(0))| \, |S(0)| \, d\xvec \right)
 + \: \frac{2}{3} \int_0^{T} \int_\Omega \phi \, 
	\frac{\partial \gamma}{\partial t} |S|^{3/2} \, d\xvec \, dt \: . $$
\end{lmm}
\begin{proof}
We prolongate $S$ to a function $\tilde{S}$, defined on $[-T,2T]$, by
$$ \tilde{S}(t) := 
   \left\{ \begin{array}{l@{\quad{\mathrm{if}}\quad}r@{\:\le t \le\:}l}
				S(-t) \: , & -T & 0 \: , \\
				S(t) \: , & 0\, & T \: , \\
				S(2T-t) \: , & T & 2T \: ,
			   \end{array} \right. $$
Owing to the corresponding properties of $S$, we can conclude that 
$\tilde{S} \in \mathrm{C} \! \left( [-T,2T]; \xLn{{3/2}}(\Omega) \right)$
and $(\partial /\partial t)\rho(\tilde{S}) 
     \in \xLinfty \! \left(-T,2T;\xWn{{-1,3}}(\Omega) \right)$.
For $\Delta t > 0$ we define
$$ X_{\Delta t} := \frac{1}{\Delta t} \int_0^{T} \int_\Omega \phi 
	\left( \rho(\tilde{S}(t)) - \rho(\tilde{S}(t-\Delta t)) \right) \,
	\tilde{S}(t) \, d\xvec \: . $$
In the limit $\Delta t \to 0$ this expression tends to
(\cf \cite[proof of Lemma~1.2]{Raviart:70})
$$ \lim_{\Delta t \to 0} X_{\Delta t} = \int_0^{T} 
   \left\langle \phi \frac{\partial \rho(S)}{\partial t} , S \right\rangle dt \: . $$
Like in the proof of Lemma~\ref{lemma:S-est:hom.Dir.pr}
an application of Young's inequality yields
\begin{eqnarray*}
X_{\Delta t} & \ge & \frac{1}{\Delta t} \int_0^{T}
	\frac{1}{3} \int_\Omega \phi \, |\rho(\tilde{S}(t))| \,
				|\tilde{S}(t)| \, d\xvec
       - \frac{1}{3} \int_\Omega \phi \, |\rho(\tilde{S}(t-\Delta t))| \,
				 |\tilde{S}(t-\Delta t)| \, d\xvec \, dt \\
& & + \, \frac{1}{\Delta t} \int_0^{T} \frac{2}{3} 
	\int_\Omega \left( \phi \, \gamma(t)
			  - \phi \, \gamma(t-\Delta t) \right)
		|\tilde{S}(t)|^{3/2} \, d\xvec \, dt \\
& = & \frac{1}{3} \frac{\phi}{\Delta t} \left( 
	\int_{T-\Delta t}^{T} \int_\Omega 
		|\rho(\tilde{S}(t))| |\tilde{S}(t)| \, d\xvec \, dt
       - \int_{-\Delta t}^0 \int_\Omega 
		|\rho(\tilde{S}(t))| |\tilde{S}(t)| \, d\xvec \, dt \right) \\
& & + \, \frac{2}{3} \int_0^{T} \int_\Omega
 	\frac{\phi}{\Delta t} \left( \gamma(t) - \gamma(t-\Delta t) \right)
		|\tilde{S}(t)|^{3/2} \, d\xvec \, dt \: .
\end{eqnarray*}
Therefore we obtain in the limit $\Delta t \to 0$
\begin{eqnarray*}
\int_0^{T} \left\langle \phi \frac{\partial \rho(S)}{\partial t} , S \right\rangle dt 
= \lim_{\Delta t \to 0} X_{\Delta t} 
& \ge & \frac{1}{3} \left( 
    \int_\Omega \phi \, |\rho(S(T))| \, |S(T)| \, d\xvec
   - \int_\Omega \phi \, |\rho(S(0))| \, |S(0)| \, d\xvec \right) \\
& & + \, \frac{2}{3} \int_0^{T} 
   \int_\Omega \phi \, \frac{\partial \gamma}{\partial t} \, |S(t)|^{3/2} \, d\xvec \, dt \: .\end{eqnarray*}
Applying the same transformations and estimations to
$$ Y_{\Delta t} := \frac{1}{\Delta t} \int_0^{T} \int_\Omega
   \phi \left( \rho(\tilde{S}(t+\Delta t)) - \rho(\tilde{S}(t)) \right) \,
	 \tilde{S}(t) \, d\xvec \: , $$
we find
\begin{eqnarray*}
\int_0^{T} \left\langle \phi \frac{\partial \rho(S)}{\partial t} , S \right\rangle dt 
= \lim_{\Delta t \to 0} Y_{\Delta t}
& \le & \frac{1}{3} \left( 
    \int_\Omega \phi |\rho(S(T))| |S(T)| \, d\xvec
   - \int_\Omega \phi |\rho(S(0))| |S(0)| \, d\xvec \right) \\
& & + \, \frac{2}{3} \int_0^{T} 
	\int_\Omega \phi \frac{\partial \gamma}{\partial t} |S(t)|^{3/2} \, d\xvec \, dt \: .
\end{eqnarray*}
Together with the estimate from the consideration of $X_{\Delta t}$ above,
the assertion follows.
\end{proof}

\begin{prpstn}					\label{prop:g=G(m)}
The limits $\mvec$ of 
$\left( \Pi_{\Delta t} \mvec_{\Delta t} \right)_{\Delta t}$
and $\gvec$ of 
$\left( \Pi_{\Delta t} G_{\Delta t} (\mvec_{\Delta t}) \right)_{\Delta t}$
satisfy $\gvec = G(\mvec)$, i.e., for all 
$\vvec \in \xLone \! \left( 0,T; (\xLn{3}(\Omega))^n \right)$:
$$
\int_0^{T} \int_\Omega \gvec \cdot \vvec \, d\xvec \, dt = 
\int_0^{T} \int_\Omega G(\mvec) \cdot \vvec \, d\xvec \, dt \: .
$$
\end{prpstn}
\begin{proof}
Again, we employ (\ref{proof:S-est:hom.Dir.pr}),
replacing $a^k(\mvec^k,\mvec^k)$ by
$\int_\Omega G^k(\mvec^k) \cdot \mvec^k \, d\xvec$, i.e.,
$$ \int_\Omega G^k(\mvec^k) \cdot \mvec^k \, d\xvec
   + \int_\Omega \frac{\phi}{\Delta t} 
	\left( \rho^k(S^k) - \rho^{k-1}(S^{k-1}) \right) S^k \, d\xvec =
   \int_\Omega f^k S^k \, d\xvec $$
and the inequality (see the proof of Lemma \ref{lemma:S-est:hom.Dir.pr})
\begin{eqnarray*}
\frac{1}{3} \int_\Omega \frac{\phi}{\Delta t} 
			|\rho^k(S^k)| |S^k| \, d\xvec
- \frac{1}{3} \int_\Omega \frac{\phi}{\Delta t}
			  |\rho^{k-1}(S^{k-1})| |S^{k-1}| \, d\xvec & & \\
+ \: \frac{2}{3} \int_\Omega \frac{\phi}{\Delta t}
		\left( \gamma^k - \gamma^{k-1} \right) |S^k|^{3/2} \, d\xvec 
& \le & \int_\Omega \frac{\phi}{\Delta t} 
	 \left( \rho^k(S^k) - \rho^{k-1}(S^{k-1}) \right) S^k \, d\xvec \: .
\end{eqnarray*}
Multiplying with $\Delta t$ and summing up for $k=1,\ldots,K$, we obtain
\begin{eqnarray*}
\frac{1}{3} \int_\Omega \phi \, |\rho^K(S^K)| |S^K| \, d\xvec
- \frac{1}{3} \int_\Omega \phi \, |\rho^0(S^0)| |S^0| \, d\xvec & & \\
+ \: \frac{2}{3} \int_{\Delta t}^{T} \int_\Omega \phi
	\left( \Pi_{\Delta t} \gamma_{\Delta t}(t)
	    - \Pi_{\Delta t} \gamma_{\Delta t}(t-\Delta t) \right)
		|\Pi_{\Delta t} S_{\Delta t}|^{3/2} \, d\xvec \, dt & & \\
+ \int_0^{T} \int_\Omega \Pi_{\Delta t} G_{\Delta t} (\mvec_{\Delta t}) \cdot
		           \Pi_{\Delta t} \mvec_{\Delta t} \, d\xvec \, dt
& \le & \int_0^{T} \int_\Omega \Pi_{\Delta t} f_{\Delta t}
			\Pi_{\Delta t} S_{\Delta t} \, d\xvec \, dt \: .
\end{eqnarray*}
Taking the limes inferior, we can conclude that
\begin{eqnarray*}
\frac{1}{3} \int_\Omega \phi \, |\rho(S(T))| |S(T)| \, d\xvec
- \frac{1}{3} \int_\Omega \phi \, |\rho(S(0))| |S(0)| \, d\xvec \\
+ \: \frac{2}{3} \int_0^{T} \int_\Omega \phi \,
	\frac{\partial \gamma}{\partial t} |S|^{3/2} \, d\xvec \, dt
+ \liminf_{\Delta t \to 0}
  \int_0^{T} \int_\Omega \Pi_{\Delta t} G_{\Delta t} (\mvec_{\Delta t}) \cdot
		       \Pi_{\Delta t} \mvec_{\Delta t} \, d\xvec \, dt
& \le & \int_0^{T} \int_\Omega f S \, d\xvec \, dt \: .
\end{eqnarray*}
Thus Lemma \ref{lemma:int(dR/dt,S)} yields the following inequality
$$ \int_0^{T} 
   \left\langle \phi \frac{\partial \rho(S)}{\partial t} , S \right\rangle dt
  + \liminf_{\Delta t \to 0}
    \int_0^{T} \int_\Omega \Pi_{\Delta t} G_{\Delta t} (\mvec_{\Delta t}) \cdot
		           \Pi_{\Delta t} \mvec_{\Delta t} \, d\xvec \, dt
  \le \int_0^{T} \int_\Omega f S \, d\xvec \, dt \: . $$
On the other hand, Proposition \ref{prop:distr.var.form} implies
$$ \int_0^{T} \left\langle \phi \frac{\partial \rho(S)}{\partial t} , S \right\rangle dt
   + \int_0^{T} \int_\Omega \gvec \cdot \mvec \, d\xvec \, dt
   = \int_0^{T} \int_\Omega f S \, d\xvec \, dt \: , $$
since from Proposition \ref{prop:S'=dS/dt} c) and d) we have
$$ \int_0^{T} \left\langle \Div(\mvec) , S \right\rangle dt =
   - \int_0^{T} \int_\Omega \mvec \cdot \nabla S \, d\xvec \, dt =
   \int_0^{T} \int_\Omega \mvec \cdot \gvec \, d\xvec \, dt \: . $$
Consequently,
$$ \liminf_{\Delta t \to 0}
   \int_0^{T} \int_\Omega \Pi_{\Delta t} G_{\Delta t} (\mvec_{\Delta t}) \cdot
		       \Pi_{\Delta t} \mvec_{\Delta t} \, d\xvec \, dt
   \le \int_0^{T} \int_\Omega \gvec \cdot \mvec \, d\xvec \, dt \: . $$
Thus we have shown that 
for arbitrary $\vvec \in \xLinfty \! \left( 0,T; (\xLn{3}(\Omega))^n \right)$
$$ 
\int_0^{T} \int_\Omega \left( \gvec - G(\vvec) \right) \cdot 
			 \left( \mvec - \vvec \right) d\xvec \, dt \ge
\liminf_{\Delta t \to 0}
   \int_0^{T} \int_\Omega \left( \Pi_{\Delta t} G_{\Delta t} (\mvec_{\Delta t})
			      - \Pi_{\Delta t} G_{\Delta t} (\vvec) \right)
			\cdot \left( \Pi_{\Delta t} \mvec_{\Delta t}
				    - \vvec \right) d\xvec \, dt \ge 0 \: .
$$ 
Now the assertion follows from the fact that $G$ is a maximal monotone 
operator on $\xLinfty \! \left( 0,T; (\xLn{3}(\Omega))^n \right)$
(\cf the proof of Thm.~1.1 in \cite{Raviart:70}).
\end{proof} 

Now we are in a position to formulate and prove our main result:
\begin{thrm}					\label{thrm:solv.cont.pr}
For all $f \in \xLinfty \left( 0,T;\xLn{3}(\Omega) \right)$
that are Lipschitz continuous in $t$ there exists a pair
$\left( \mvec , S \right) \in \xLinfty \! \left( 0,T; (\xLn{3}(\Omega))^n \right)
 \times \xLinfty \big( 0,T; W_0^{1,3/2}(\Omega)) \big)$ such that
$$ \begin{array}{rcl@{\quad\mbox{for all}\quad}l}
\displaystyle \int_0^{T} \int_\Omega G(\mvec) \cdot \vvec \, d\xvec \, dt
- \int_0^{T} \left\langle \Div(\vvec) , S \right\rangle dt & = & 0 
& \vvec \in \xLone \! \left( 0,T; (\xLn{3}(\Omega))^n \right) \: , \\[2.5ex]
\displaystyle \int_0^{T} \left\langle \phi \, \frac{\partial \rho(S)}{\partial t} , q \right\rangle dt
 + \int_0^{T} \left\langle \Div(\mvec) , q \right\rangle dt
& = & \displaystyle \int_0^{T} \int_\Omega f \, q \, d\xvec \, dt
& q \in \xLone \big( 0,T; W_0^{1,3/2}(\Omega)) \big) \; .
\end{array} $$
\end{thrm}
\begin{proof}
Let $\mvec$ be the limit of 
$\left( \Pi_{\Delta t} \mvec_{\Delta t} \right)_{\Delta t}$ and 
$S$ be the limit of $\left( \Pi_{\Delta t} S_{\Delta t} \right)_{\Delta t}$.
Then Proposition \ref{prop:g=G(m)} and Proposition \ref{prop:S'=dS/dt} d) 
imply that
$$ 
\int_0^{T} \int_\Omega G(\mvec) \cdot \vvec \, d\xvec \, dt =
\int_0^{T} \int_\Omega \gvec \cdot \vvec \, d\xvec \, dt =
- \int_0^{T} \int_\Omega \nabla S \cdot \vvec \, d\xvec \, dt =
\int_0^{T} \left\langle \Div(\vvec) , S \right\rangle dt
$$ 
for all $\vvec \in \xLone \left( 0,T; (\xLn{3}(\Omega))^n \right)$.
Thus $(\mvec,S)$ satisfies the first equation above.
In Proposition \ref{prop:distr.var.form} we have seen that $(\mvec,S)$
fulfills the second equation, too.
\end{proof}

\begin{rmrk}
\begin{enumerate} \renewcommand{\labelenumi}{\alph{enumi})}
\item By means of the definition of the generalized divergence,
    we can replace the dual pairing $\langle \Div(\vvec) , q \rangle$
    for $\vvec \in \left( \xLn{3}(\Omega) \right)^n$ and
    $q \in \xWn{{1,3/2}}_0(\Omega)$ with the integral 
    $-\int_\Omega \vvec \cdot \nabla q \, d\xvec$.
    Thus, in the case of the continuous transient problem,
    we have established the existence of a solution 
    of the primal mixed formulation (\cf \cite[Sect.~I.3.2]{Roberts/Thomas}).
    In contrast, we considered the uniqueness and existence of a solution
    of the dual mixed formulation for the stationary
    and semi-discrete transient Problem.
    This lack of regularity of the vector solution $\mvec$ hinders
    the consideration of more general boundary conditions.
\item Assuming additional regularity properties of the solution,
    Amirat \cite{Amirat:91} showed that the solution
    to the corresponding parabolic Neumann-problem is unique.
    Furthermore, he proved that the solution is positive
    provided that the initial and boundary conditions satisfy 
    corresponding requirements.    
\end{enumerate}
\end{rmrk}

\begin{appendix}
\section{Properties of $\xWn{s}(\Div;\Omega)$}
\label{app:W^s(div)}
We introduce the generalization $\xWn{s}(\Div;\Omega)$ of 
$\mathrm{H}(\Div;\Omega)$, defined by
$$ \xWn{s}(\Div;\Omega) := 
   \left\{ \vvec \in \left( \xLn{s}(\Omega) \right)^n \bigm|
   	   \Div(\vvec) \in \xLn{s}(\Omega) \right\} \: , $$
and equip it with the norm
$$ \|\vvec\|_{\xWn{s}(\Div;\Omega)} := 
   \left( \int_\Omega \sum_{i=1}^n |v_i(\xvec)|^s \, d\xvec + 
          \int_\Omega |\Div(\vvec(\xvec))|^s \, d\xvec \right)^{1/s} \ , $$
where $\vvec = (v_1,\ldots,v_n)^T$.
Since $\xWn{s}(\Div;\Omega)$ is a closed subspace 
of $\left( L^s(\Omega) \right)^{n+1}$, it follows that
$\xWn{s}(\Div;\Omega)$ is a reflexive Banach space.

It is straightforward to extend the proofs of Thm.~2.4 and Thm.~2.5 
in \cite{Girault/Raviart} to show the next two lemmas:
\begin{lmm}					\label{dicht-lemma}
The space $\mathcal{D}(\bar{\Omega})^n$ is dense in $W^s(\Div;\Omega)$.
\end{lmm}

\begin{lmm}
The mapping $\gamma_n : \vvec \mapsto \vvec \cdot \nvec$
defined on $\mathcal{D}(\bar{\Omega})^n$ can be extended by continuity
to a linear and continuous mapping, still denoted by $\gamma_n$,
from $W^s(\Div;\Omega)$ into $\left( W^{1/s,r}(\partial \Omega) \right)'$.
In particular, Green's formula
\begin{equation}			\label{Green's}
   \int_\Omega \vvec \cdot \nabla \psi \, d\xvec +
   \int_\Omega \Div(\vvec) \psi \, d\xvec = 
   \int_{\partial\Omega} \psi (\vvec \cdot \nvec) \, d\sigma
\end{equation}
holds for every $\vvec \in \xWn{s}(\Div;\Omega)$ and 
$\psi \in \xWn{{1,r}}(\Omega)\,$, where $1/s + 1/r = 1$.
\end{lmm}

For $s>1$, the well known inf-sup condition 
(see e.g.\ \cite[{\S}II.1]{Brezzi/Fortin}) can be extended, too. 
Generalizing the definition of the bilinear form $b$ 
from Section \ref{sec:stat.pr} onto 
$\xWn{s}(\Div;\Omega) \times \xLn{r}(\Omega)$, we define
$b(\vvec,q) := \int_\Omega \Div(\vvec) \, q \, d\xvec$
for $\vvec \in \xWn{s}(\Div;\Omega)$, $q \in \xLn{{r}}(\Omega)$.

\begin{lmm}
Let $s>1$ and $1/s+1/r=1$. 
Then there exists a constant $\theta > 0$ such that
\begin{equation}				\label{inf-sup}
\theta \|q\|_{0,r,\Omega} \le \sup_{\vvec \in \xWn{s}(\Div;\Omega)} 
\frac{b(\vvec,q)}{\|\vvec\|_{\xWn{s}(\Div;\Omega)}} 
\quad \mbox{for all } \: \vvec \in \xWn{s}(\Div;\Omega) \, , ~
		q \in \xLn{{r}}(\Omega)\: .
\end{equation}
\end{lmm}
\begin{proof}
We define a mapping $B : \xWn{s}(\Div;\Omega) \to 
\xLn{s}(\Omega) = \left( \xLn{{r}}(\Omega) \right)'$ by means of
$\left\langle B \vvec , q \right\rangle = b(\vvec,q)$.
Since $s>1$, for every $p \in \xLn{s}(\Omega)$ there exists 
a Newtonian potential $N_p \in \xWn{{2,s}}(\Omega)$ such that $\Delta N_p = p$
almost everywhere. Then $\vvec := \nabla N_p$ satisfies
$\vvec \in \xWn{s}(\Div; \Omega)$
and $\Div(\vvec) = \Delta N_p = p$ in $\xLn{s}(\Omega)\,$.
Therefore $B$ is onto.
Applying Lemma A.1 of \cite{Sandri:93}, the assertion follows.
\end{proof}
\end{appendix}

%

\end{document}